\documentclass[11pt]{article}

\usepackage[dvips]{graphicx}
\usepackage{epsfig,amsmath,latexsym}
\usepackage{amsfonts}

\newtheorem{theorem}{Theorem}[section]

\newtheorem{defi}{Definition}[section]
\newtheorem{lemma}[theorem]{Lemma}

\catcode`\@=11
\renewcommand{\section}{
        \setcounter{equation}{0}
        \@startsection {section}{1}{\z@}{-3.5ex plus -1ex minus
        -.2ex}{2.3ex plus .2ex}{\large\bf}
        }
\catcode`@=12

\def\slfrac#1#2{\hbox{\kern.1em %
 \raise.5ex\hbox{\the\scriptfont0 #1}\kern-.11em %
 /\kern-.15em\lower.25ex\hbox{\the\scriptfont0 #2}}}

\newcommand{\Aut}{{\operatorname{Aut}}}

\newcommand{\beq}{\begin{eqnarray}}
\newcommand{\eeq}{\end{eqnarray}}
\newcommand{\beql}[1]{\begin{eqnarray}\label{#1}}
\newcommand{\beqs}{\begin{eqnarray*}}
\newcommand{\eeqs}{\end{eqnarray*}}
\newcommand{\eqn}[1]{(\ref{#1})}

\newcommand{\rr}{{\mathbb R}}
\newcommand{\zz}{{\mathbb Z}}
\newcommand{\qq}{{\mathbb Q}}

\newcommand{\ba}{{\mathbf a}}
\newcommand{\bb}{{\mathbf b}}

\newcommand{\bv}{{\mathbf v}}

\newcommand{\bx}{{\mathbf x}}

\newcommand{\bI}{{\mathbf I}}
\newcommand{\bJ}{{\mathbf J}}
\newcommand{\bP}{{\mathbf P}}
\newcommand{\bQ}{{\mathbf Q}}
\newcommand{\bS}{{\mathbf S}}
\newcommand{\bU}{{\mathbf U}}

\newcommand{\sA}{{\mathcal A}}

\newcommand{\sD}{{\mathcal D}}

\newcommand{\sL}{{\mathcal L}}

\newcommand{\sO}{{\mathcal O}}
\makeatletter
\def\section{\@startsection {section}{1}{\z@}{-3.5ex plus -1ex minus 
 -.2ex}{2.3ex plus .2ex}{\normalsize\bf}}
\def\subsection{\@startsection {subsection}{1}{\z@}{-3.5ex plus -1ex minus
 -.2ex}{2.3ex plus .2ex}{\normalsize\bf}}
\def\@sect#1#2#3#4#5#6[#7]#8{\ifnum #2>\c@secnumdepth
     \def\@svsec{}\else
     \refstepcounter{#1}\edef\@svsec{\csname the#1\endcsname.\hskip .75em }\fi
     \@tempskipa #5\relax
      \ifdim \@tempskipa>\z@
        \begingroup #6\relax
          \@hangfrom{\hskip #3\relax\@svsec}{\interlinepenalty \@M #8\par}%
        \endgroup
       \csname #1mark\endcsname{#7}\addcontentsline
         {toc}{#1}{\ifnum #2>\c@secnumdepth \else
                      \protect\numberline{\csname the#1\endcsname}\fi
                    #7}\else
        \def\@svsechd{#6\hskip #3\@svsec #8\csname #1mark\endcsname
                      {#7}\addcontentsline
                           {toc}{#1}{\ifnum #2>\c@secnumdepth \else
                             \protect\numberline{\csname the#1\endcsname}\fi
                       #7}}\fi
     \@xsect{#5}}
\def\@begintheorem#1#2{\it \trivlist \item[\hskip \labelsep{\bf #1\ #2.}]}
\makeatother

\begin{document}
  
  \begin{center}
    {\Large {\bf Apollonian Circle Packings: Number Theory}} \\
    {\Large {\bf II. Spherical and Hyperbolic Packings}}\\
    \vspace{1.5\baselineskip}
           {\em Nicholas Eriksson} \\
           \vspace*{.2\baselineskip}
           University of California at Berkeley \\
           Berkeley, CA 94720 \\
           
           \vspace*{.2\baselineskip}
           
                   {\em Jeffrey C. Lagarias} \\
                   \vspace*{.2\baselineskip}
                   University of Michigan\\
                   Ann Arbor, MI 48109\\
                   \vspace*{1.5\baselineskip}
                  
                   (February 1, 2005)
                   
                   \vspace{1.5\baselineskip}
                          {\bf ABSTRACT}
  \end{center}
  
  Apollonian circle packings arise by repeatedly filling the interstices
  between mutually tangent circles with further tangent circles.  In
  Euclidean space it is possible for every circle in such a packing to
  have integer radius of curvature, and we call such a packing an {\em
    integral Apollonian circle packing.} There are infinitely many
  different integral packings; these were studied in the paper
  \cite{GLMWY21}.  Integral circle packings also exist
  in spherical and hyperbolic space, provided a suitable definition of
  curvature is used (see \cite{LMW02}) and again there are an infinite
  number of different integral packings. This paper studies
  number-theoretic properties of such packings. This amounts to 
  studying  the orbits of a particular subgroup $\sA$ of the
  group of integral automorphs of the indefinite quaternary quadratic form
  $Q_{\sD}(w, x, y, z)= 2(w^2+x^2 +y^2 + z^2) - (w+x+y+z)^2$.
  This subgroup, called the Apollonian group, 
  acts on integer solutions $Q_{\sD}(w, x, y, z)=k$.
  This paper  gives a reduction theory for orbits of $\sA$ 
  acting on integer solutions to $Q_{\sD}(w, x, y, z)=k$  
  valid for all integer $k$. It also classifies orbits for all
  $k \equiv 0 \pmod{4}$ in terms of an extra parameter $n$ and
  an auxiliary class group (depending on $n$ and $k$), and 
  studies congruence conditions on integers in a given orbit.
  
  \vspace*{1.5\baselineskip}
  

%
  
\setlength{\baselineskip}{1.0\baselineskip}

%
%
%
%
\section{Introduction}
\begin{figure}
\begin{center}
\includegraphics[height=3.5in  ]{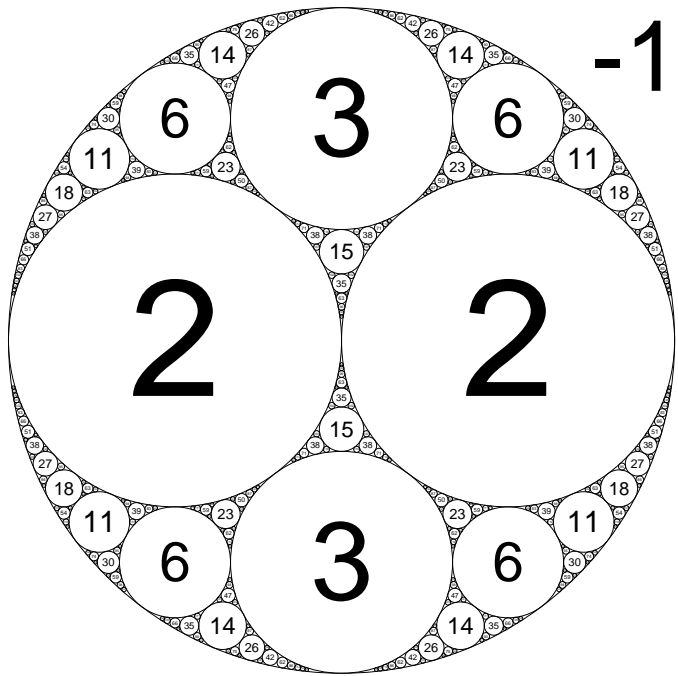}
\end{center}
\caption{The Euclidean Apollonian circle packing ($-1$,$2$,$2$,$3$)}
\label{fig0}
\end{figure}

\setcounter{equation}{0}
A Descartes configuration is a set of four mutually touching
circles with distinct tangents.  
A (Euclidean) Apollonian circle packing in the plane is constructed starting
from a Descartes configuration  by recursively
adding circles tangent to three of the circles already constructed in the
packing; see \cite{LMW02}, \cite{GLMWY11}, \cite{GLMWY12}, 
and Figure~\ref{fig0}. 

Apollonian packings can be  completely described by
the set of Descartes configurations they contain.
The curvatures $(a,~b,~c,~d)$ of the circles in a
Descartes configuration satisfy 
an algebraic relation called the Descartes circle theorem.
Using the 
{\em Descartes quadratic form}
\begin{equation}~\label{eq0}
Q_{\sD}(a, b, c, d) := 2(a^2+b^2+c^2+d^2)-(a+b+c+d)^2, 
\end{equation}
the Descartes circle theorem
states that
\beq\label{descartes}
Q_{\sD}(a,b,c,d) = 0.
\eeq
Descartes originally stated this relation in another 
algebraic form; see \cite{LMW02}.
In order to have this formula work in all cases, the curvatures 
of the circles in a
Descartes configuration must be given proper signs, as specified
in \cite{LMW02} and \cite{GLMWY11}.

The form $Q_{\sD}$ is  an integral quadratic form in four variables,
and the equation \eqn{descartes} has many integer solutions.
If the initial Descartes configuration generating a (Euclidean)
Apollonian packing has all four circles
with integral curvatures, then it turns out that 
all the circles in the packing have integer curvatures.
This integer structure was studied  in the
first paper in this series \cite{GLMWY21}
using the action of a particular
discrete subgroup of automorphisms of the Descartes quadratic
form, which was there termed  the {\em Apollonian group} $\sA$.
It is a group of integer matrices, with
generators explicitly given in \S3. 
The Descartes configurations in  a (Euclidean)  Apollonian packing 
are described by an
orbit of the Apollonian group.
That is, the Apollonian group acts on the column vector of curvatures
of the initial Descartes configuration, and  generates the
curvature vectors of all Descartes configurations in the 
packing. This fact explains the preservation of integrality 
of the curvatures.  The
Apollonian group and its relevance for Apollonian packings was already
noted in 1992 by S\"{o}derberg \cite{So92}. 

The  first paper in this
series \cite{GLMWY21} studied  Diophantine properties of the
curvatures in Euclidean integer Apollonian circle
packings, based on the Apollonian group action.
In \cite{LMW02} it was observed that 
there exist notions of 
integer Apollonian packings in spherical and
hyperbolic geometry. These are based on  analogues 
of the Descartes equation valid in these geometries,
as follows.
In hyperbolic space, if the ``curvature'' of a 
circle with hyperbolic radius $r$ is
taken to be $\coth(r)$, then the modified Descartes equation 
for the ``curvatures'' of circles in a Descartes
configuration is 
\beq\label{hypdescartes}
Q_{\sD}(a, b, c, d)=  4.
\eeq
Similarly, in spherical space, if the ``curvature'' of a circle with 
spherical radius
$r$ is defined  to be $\cot(r)$, then the 
modified Descartes equation is
\beq\label{sphdescartes} 
Q_{\sD}(a, b, c, d)=  - 4.  
\eeq
There exist integral Apollonian packings in both these
geometries, in the sense that 
 the  (spherical or hyperbolic) ``curvatures'' of
all circles in the packing are integral; see  the examples
in Figures~\ref{fig:spherical} and \ref{fig:hyper} below.
Again, the  vectors of ``curvatures'' of all Descartes configurations
in such a packing form an orbit under the action
of the Apollonian group.

The purpose of this paper is to study
Diophantine properties of
the orbits  of the Apollonian group $\sA$ acting 
on the integer solutions to the Diophantine 
equation
\begin{equation}~\label{103a}
Q_{\sD}(a, b, c, d)= k,
\end{equation}
where $k$ is  a fixed integer. Integer solutions exist for
  $k \equiv 0$ or $3~ (\bmod~4)$, and our  results on
reduction theory in \S3 are proved in this generality.
Our other results are established under the 
more restrictive condition 
$k \equiv 0 \pmod{4}$, where we
will always write $k=4m$. These results cover the 
motivating cases $k= \pm 4$
corresponding to spherical and hyperbolic Apollonian packings,
as well as the Euclidean case $k=0$. 
We note that the  Apollonian group $\sA$ is a subgroup of 
{\em infinite} index in the
group $\Aut(Q_{\sD}, \zz)$ of all integer automorphs of
the Descartes quadratic form $Q_{\sD}$. As a consequence, the
integral solutions of $Q_{\sD}(a, b, c, d)= k$ (for
$k \equiv 0$ or $3~ (\bmod~4)$) fall in  an infinite
number of orbits under the  action of $\sA$. In contrast, these same 
integral points fall in a finite number of orbits under the  full 
$\Aut(Q_{\sD}, \zz)$-action.

In \S \ref{integral}, we
study the distribution
of  integer solutions to the equation $Q_{\sD}(a,b,c,d) = 4m$,
determining  asymptotics for the number
of such  solutions to these equations of size below a given bound.
These asymptotics are obtained via a bijective correspondence between 
integral representations of integers $4m$ by the
Descartes form and representations of integers $2m$ by the Lorentzian form
$Q_{\sL}(W, X, Y, Z) = -W^2 + X^2 + Y^2+ Z^2$, and then applying
results  of Ratcliffe and Tschantz \cite{RT97}
on integer solutions to the Lorentzian form.

\begin{figure}
\begin{center}
\includegraphics[height=3.5in  ]{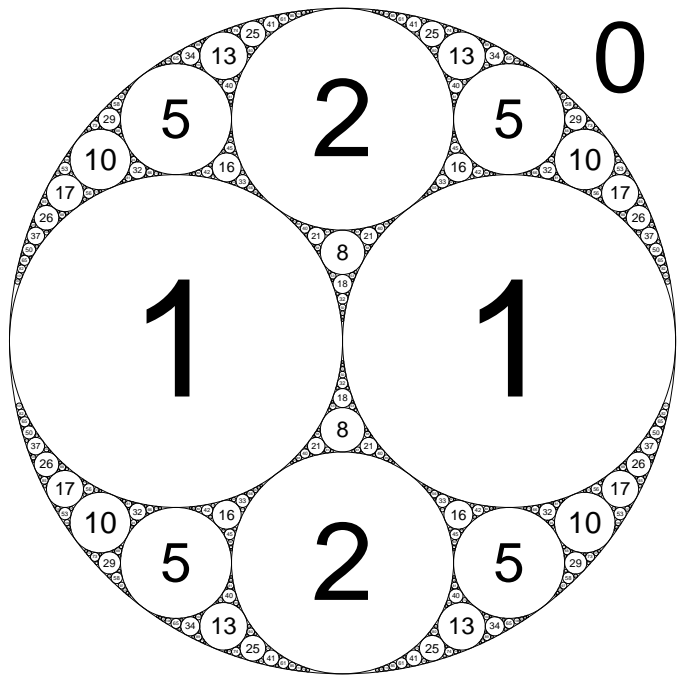}
\end{center}
\caption{The spherical Apollonian circle packing ($0$,$1$,$1$,$2$)}
\label{fig:spherical}
\end{figure}

\begin{figure}
\begin{center}
\includegraphics[height=3.5in  ]{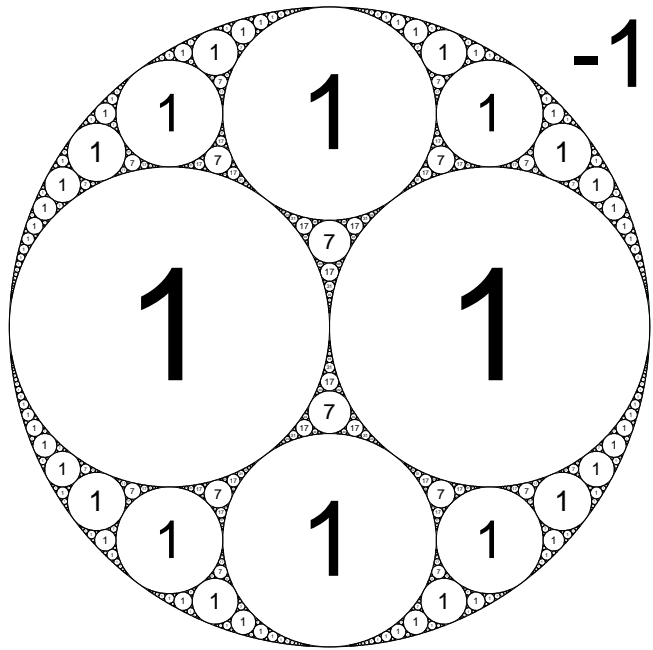}
\end{center}
\caption{The hyperbolic Apollonian circle packing ($-1$,$1$,$1$,$1$)}
\label{fig:hyper}
\end{figure}

In \S \ref{reduction}, we describe the Apollonian group,
and apply it to give a reduction algorithm for
quadruples $\bv=(a, b, c,d)$ satisfying
$Q_{\sD}(\bv) = k$ in a fixed orbit of 
the Apollonian group, starting from an
initial integer quadruple satisfying  $Q_{\sD}(\bv) = k$.
This reduction algorithm differs
in appearance  from that given in part I \cite[Sect. 3]{GLMWY21}, which
treated the case $k=0$; however we show its steps coincide 
with those of the earlier algorithm when $k=0$. 
A quadruple obtained at the
end of this algorithm is called a {\em reduced quadruple}. 
We classify
reduced quadruples as two types, 
{\em root quadruples} and {\em exceptional quadruples}.
We show that for each orbit containing a root quadruple, 
the root quadruple  is
the unique reduced quadruple in that orbit.
For each fixed $k$ we show
that there are at most a  finite number of orbits
not containing a  root quadruple.
Each of these orbits contain at least one and
at most finitely many reduced quadruples,
each of which is by definition an exceptional quadruple.
We show that such orbits can only 
occur when $k> 0$, and we  call them {\em exceptional orbits}.
Their existence is  an  interesting new phenomenon 
uncovered in this work. They do occur in the
hyperbolic Apollonian packing case  $k=4$,
where we show there are exactly two exceptional orbits, labeled by  the
hyperbolic Descartes  quadruples
$(0,0,0, 2)$ and $(-1,0,0,1)$, respectively. We give a geometric
interpretation of these exceptional 
hyperbolic Apollonian packings.

As examples, 
Figure~\ref{fig:spherical} shows the
spherical ($k = -4$) Apollonian packing with root quadruple $(0,1,1,2)$, 
after stereographic projection.
Figure~\ref{fig:hyper} shows
the hyperbolic ($k=4$) Apollonian packing with root quadruple $(-1,1,1,1),$
using the unit disc model of the hyperbolic plane.
For more on the geometry of
Apollonian packings in spherical and hyperbolic space 
consult \cite{LMW02}.

In \S \ref{root}, we count the number $N_{\text{root}}(4m, -n)$ of 
root quadruples $\ba= (a, b,c,d)$
with $Q_{\sD}(\ba) = 4m$ 
having  smallest member  $a= -n$. We give a formula
for this number in terms of the 
class number of (not necessarily primitive)
integral binary quadratic forms of discriminant $-4(n^2 - m)$ 
under $GL(2,\zz)$-equivalence.
The class number 
interpretation applies only when $n^2 > m$.
For certain $m>0$ there  may exist root quadruples with
parameters $n^2 \le m$, and information on them
is given in Theorem~\ref{th33} in \S \ref{reduction}.
For the hyperbolic case the only value not covered
is  $n=-1$.  Associated
to it is an infinite family of distinct root quadruples
$\{(-1, 1, c, c)|~ c \ge 1\}$. Figure~\ref{fig:hyper} pictures
the Apollonian packing with root quadruple 
$(-1,1,1,1)$.

The class number equivalence
 allows us to derive good upper  bounds  for the number of
such root quadruples with smallest value $a=-n$ 
as $n \to \infty$, namely
\[
N_{\text{root}}(4m, -n ) = O \left( n (\log n) (\log\log n)^2 \right). 
\]
If $m= \pm 1$, then  these  bounds are interpretable
as bounds on the 
integral spherical and hyperbolic packings with boundary an
enclosing  circle of (spherical or hyperbolic) curvature $n$.

In \S \ref{congruences} we investigate congruence restrictions
on  the integer
curvatures which occur in a packing, for the case of spherical
and hyperbolic packings, i.e., $k = \pm 4$. We show that there 
are always non-trivial 
congruence conditions modulo $12$
on the set of integers that can occur in a given packing.
It seems reasonable to expect
 that such congruence conditions can only involve powers of
the primes 2 and 3 and that all sufficiently large integers 
which are not excluded by 
these congruence conditions actually occur.

Various open problems are  raised at the end of \S \ref{root}
and in \S \ref{congruences}. In addition, in the case of 
Euclidean Apollonian
packings, it is known that there are Apollonian packings which have
stronger integer properties, involving the centers of
the circles as well as the curvatures. These were studied 
in  Graham et al.\ \cite{GLMWY11}, \cite{GLMWY12}, and \cite{GLMWY13}. 
It remains to be seen
if there are analogues of such properties (for circle centers)
in the spherical and hyperbolic cases.

\paragraph{Acknowledgments.}
Much of this work was done while the authors
were at AT\&T Labs-Research, 
whom the authors thank for support.  N.~Eriksson was
also supported by an NDSEG fellowship.
The authors thank the reviewer for helpful comments.

%
%
%
%
%
\section{Integral Descartes Quadruples}
\label{integral}
\setcounter{equation}{0}

The {\em Descartes quadratic form} $Q_{\sD}$ is the 
quaternary quadratic form
\beq \label{desform}
Q_\sD(w,x,y,z) = 2(w^2 + x^2 + y^2 + z^2) - (w + x + y + z)^2.
\eeq
Its matrix representation, for $\bv= [w, x, y, z]^T$, is
\beqs
Q_\sD( \mathbf v) = \mathbf v^T \bQ_\sD \mathbf v = \mathbf v^T 
\begin{pmatrix} 
1& -1& -1& -1\\
 -1& 1& -1& -1\\
 -1& -1& 1& -1\\
 -1& -1& -1& 1
\end{pmatrix}
\mathbf v.
\eeqs
This form is indefinite with signature $(+,+,+,-)$.
We consider integral representations of an integer $k$ by
this quaternary quadratic form. 
For an integer quadruple  $\bv=(w, x, y, z) \in \zz^4$,
$Q_\sD(w,x,y,z) \equiv 0$ or $3~(\bmod~4)$, according to the parity
of $w+x+y+z$. 
We define the {\em Euclidean height} $H(\mathbf v)$
of a (real) quadruple $\bv= (w, x, y,z)$ to be  
\beqs
H(\mathbf v) := (w^2 + x^2 + y^2 + z^2)^{1/2}.
\eeqs
Notice that we reserve the usual notation, $|v|$, for another meaning in 
\S \ref{reduction}.

We relate integer representations of
even $k$ of the Descartes form to integer
representations of $\frac{k}{2}$ of
the {\em Lorentz form}
\beqs
Q_\sL (W, X, Y, Z ) = -W^2 + X^2 + Y^2 + Z^2.
\eeqs
in Lemma~\ref{lem:lorentz} below.
The relation only works  for even $k$, 
so in this section we consider only the 
case $k\equiv 0 \pmod{4}$, so that
\begin{equation}~\label{repn}
Q_\sD(w,x,y,z) = k = 4m.
\end{equation}
The Lorentz form has
the  matrix representation, for $\mathbf v = [W,X,Y,Z]^T$,
\beqs
Q_\sL(\mathbf v) =\mathbf v^T \begin{pmatrix}
-1& 0 & 0 & 0\\
0& 1 & 0 & 0\\
0& 0 & 1 & 0\\
0& 0 & 0 & 1
\end{pmatrix}
\mathbf v.
\eeqs
The Descartes form and
Lorentz form are related by 
\beq~\label{deslor}
2Q_\sL = \bJ_0 Q_\sD \bJ_0^T, 
\eeq
where
\beqs
 \bJ_0  = \frac 1 2 \left(
\begin{array}{rrrr}
1& 1 &  1 & 1\\
1& 1 & -1 & -1\\
1& -1 & 1 & -1\\
1& -1 & -1 & 1\\
\end{array}
\right).
\eeqs
Notice that $\bJ_0^2 = \bI$.  

Let $N_\sL(k,T)$ count the number of integer representations of $k$ by
the Lorentz form of Euclidean height at most $T$ and
let $N_{\sD}(k, T)$ denote the number of solutions of
\eqn{repn} of Euclidean height at most $T$.

\begin{lemma}\label{lem:lorentz}
The mapping $(W,X,Y,Z)^T = \bJ_0 (w,x,y,z)^T$ gives a bijection between
real solutions of $Q_\sD (w,x,y,z) = 4r$ 
and real solutions of $Q_\sL (W,X,Y,Z) = 2r$ 
for each $r \in \rr$, and  this bijection preserves the Euclidean height
of solutions. 
In the case where $r=m \in \zz$ this bijection restricts to a
bijection between integer representations of $4m$ by the Descartes form and
integer representations of $2m$ by the Lorentzian form.
In particular, for all integers $m$, 
\beqs 
N_{\sD}(4m, T)  = N_\sL(2m,T).
\eeqs
\end{lemma}

\paragraph{Proof.}
There is a bijection on the level of real numbers,
by \eqn{deslor}, since $\bJ_0$ is
invertible, and  $Q_\sD(w,x,y,z) = 4r$ if and only if
$Q_\sL(W,X,Y,Z) = 2r$.   
A  calculation  using $\bJ_0^2 = \bI$ shows that
$$
H(w,x,y,z)= w^2+x^2 +y^2 +x^2 = W^2 + X^2 +Y^2 + Z^2 = H(W, X, Y, Z).
$$

Now suppose that $r=m \in \zz$. The mapping
takes integral solutions $(w,x,y,z)$ to integral
solutions $(W,X,Y, Z)$ because all 
integral solutions of $Q_\sD(w,x,y,z) = 4m$ satisfy $w + x + y + z
\equiv 0 \pmod 2$, as can be seen by reducing the 
Descartes equation \eqn{repn}
modulo 2. In the reverse direction integral solutions go to
integral solutions because all solutions of  $Q_\sL(w,x,y,z) = 2m$ 
have $W+X+Y+Z \equiv 0 \pmod 2$, similarly.
~~~$\Box$\\

Representation of integers by the Lorentz form has been much studied;
see Radcliffe and Tschantz \cite{RT97}. Their results
immediately yield the following result.

\begin{theorem}~\label{th21}
For each nonzero integer $m$, there is an explicitly computable
rational number $c(2m)$ such that 
\beq~\label{eq201}
N_{\sD}(4m, T) = c(2m) \frac {\pi} {L(2,\chi_{-4})}T^2 + o(T^{2}),
\eeq
in which 
\beqs
L(2,\chi_{-4}) = \sum_{n=0}^\infty \frac {(-1)^n} {(2n + 1)^2} 
\approx 0.9159.
\eeqs
We have $c(-2)= \frac{3}{4}$ and $c(2)= \frac{5}{4}.$
\end{theorem}

\begin{pf} From 
Lemma~\ref{lem:lorentz} the problem of counting integral solutions to 
$Q_{\sD}(w,x,y,z) = 4m$
of height at most $T$ is equivalent 
to counting 
the number of integral solutions to $Q_\sL(W,X,Y,Z) = 2m$ of 
Euclidean height at
most $T$.  
We use asymptotic
formulae of Radcliffe and Tschantz 
for the number $r(n,k,T)$ of solutions to the Lorentz form 
$-X_0^2 + X_1^2 + \dots + X_n^2=k$ with  $|X_0| \le T$, for
any $n \ge 2$ 
for any fixed nonzero $k$. We take $n=3$, and their
formula for $k <0$ in 
\cite[Theorem 3, p.\ 505]{RT97} is
\beq~\label{eq202}
r(3,k, T) = \frac{ \mbox{Vol}(S^2)}{2} \delta(3, k) T^2 + O( T^{3/2}),
\eeq
in which  $  \mbox{Vol}(S^2)= \pi$ and 
$\delta(3, k)$ is a density constant for the representation
in the sense of Siegel \cite{Si35}.
We need asymptotics for the number of solutions $s(n,k,T)$ to
the Lorentz form $Q_\sL(W,X,Y,Z) =k$ of 
Euclidean height below $T$, which is given by 
\beqs
s(n, k, T) := r \left(n, k, \sqrt{\frac{1}{2}(T^2 - k)}\right).
\eeqs 
The asymptotic formula  \eqn{eq202} yields  
\beq~\label{eq203} 
s(3,k, T) = \frac{ \mbox{Vol}(S^2)}{4} \delta(3, 2m) T^2 + O( T^{3/2}).
\eeq
Their result \cite[Theorem 12, p.\ 518]{RT97}
evaluates the constant $\delta(3, k)$
in terms of the value of an $L$-function, as
\beqs
\delta(3, k) = \prod_{p~|~k} 
\left(\frac{\delta_p(3, k)}{(1 - (\frac{-4}{p}) p^{-2})} \right)\cdot
\frac{1}{L(2, \chi_{-4})},
\eeqs
in which each $\delta_p(3, k)$ is an effectively 
computable  local density which is a rational
number. We now set  $k=2m$,  and  obtain a formula of
the desired shape \eqn{eq201}, with 
\beqs
c(2m) = \frac{1}{4} 
\left(\frac{\delta_p(3, 2m)}{(1 - (\frac{-4}{p}) p^{-2})} \right).   
\eeqs

Radcliffe and Tschantz \cite[Theorem 4, p.\ 510]{RT97}
also obtained an asymptotic formula
for solutions to the Lorentz form with  $k>0$, which gives
\beqs
r(3, k, T) =  \frac{ \mbox{Vol}(S^2)}{2} \delta(3, k) T^2 + o(T^2).
\eeqs
We derive a similar  asymptotic formula for $s(n,k, T)$ (without 
explicit error term), by setting  $k=2m$, and proceeding as above.

The particular values  $c(-2)=\frac{3}{4}$ and $c(2)=\frac{5}{4}$
are derived from Table II in Radcliffe and
Tschantz \cite[p. 521]{RT97}, using
$k= \pm2$.
~~~$\Box$
\end{pf}

\paragraph{Remarks.} (1)  The case $k=0$ 
exhibits a significant difference  from the pattern
of Theorem~\ref{th21}. In 
\cite[Theorem 2.1]{GLMWY21} it was shown  that
\beqs
N_{\sD}(0, T) = \frac{1}{4}  \frac {\pi^2} {L(2,\chi_{-4})}T^2 + 
O(T (\log T)^2).
\eeqs
That is, the coefficient $c(0) = \frac{\pi}{4}$ is {\em transcendental}.

(2) Theorem~\ref{th21}  implies that the number of hyperbolic  Descartes
quadruples of Euclidean height below $T$ is asymptotically
$\frac{5}{3}$ the number of spherical  Descartes quadruples of
height below $T$, as $T \to \infty.$ That is,
\beq~\label{ratiobound}
\lim_{T \to \infty} \frac{N_{\sD}(4, T)}{N_{\sD}(-4, T)} = \frac{5}{3}. 
\eeq

%
%
%
%
%
\section{Reduction Theory and Root Quadruples}
\label{reduction}
\setcounter{equation}{0}

The results of this section apply to  the general case
\beq~\label{grepn}
Q_\sD(w,x,y,z) = k,
\eeq
with $k \equiv 0$ or $3 ~(\bmod~4)$ an integer.
We introduce the Apollonian group,
and present a reduction procedure which, given an integral
Descartes quadruple, 
uses the action of the Apollonian group
to transform it to a reduced quadruple, which is
minimal according to a certain measure.
The reduced quadruples then classify integral Apollonian
circle packings, generally uniquely,
but up to a finite ambiguity in some exceptional
cases.

%
%
%
%

\subsection{Apollonian Group}

In \cite{GLMWY21} and \cite{LMW02} Apollonian circle 
packings were specified in terms of  the set of Descartes configurations
they contain. This set of configurations completely describes
the packing, and it was  observed there that, in 
a suitable coordinate system, they form 
as the orbit of any single Descartes configuration in the
packing under the motion of a discrete group,
the Apollonian group. 

\begin{defi}~\label{de31}
{\em 
The {\em Apollonian group} $\sA = \langle  \bS_1, \bS_2, \bS_3, \bS_4 \rangle$
 is the subgroup of $GL(4, \zz)$ generated by the four
integer $4 \times 4$ matrices}
\begin{gather*}
\bS_1 = 
\begin{pmatrix}
-1 & 2 & 2 & 2\\
0  & 1 & 0 & 0\\
0  & 0 & 1 & 0\\
0  & 0 & 0 & 1
\end{pmatrix},~~~
\bS_2 = 
\begin{pmatrix}
 1 & 0 & 0 & 0\\
2  &-1 & 2 & 2\\
0  & 0 & 1 & 0\\
0  & 0 & 0 & 1
\end{pmatrix},\\
\bS_3  = 
\begin{pmatrix}
 1 & 0 & 0 & 0\\
0  & 1 & 0 & 0\\
2  & 2 &-1 & 2\\
0  & 0 & 0 & 1
\end{pmatrix},~~~
\bS_4 = 
\begin{pmatrix}
 1 & 0 & 0 & 0\\
0  & 1 & 0 & 0\\
0  & 0 & 1 & 0\\
2  & 2 & 2 &-1
\end{pmatrix}.
\end{gather*}
\end{defi}

The Apollonian group is a subgroup of the  
integer automorphism group $\Aut(Q_{\sD},\zz)$ 
of the Descartes form under congruence, which is  the largest 
subgroup of $GL(4, \zz)$ which leaves $Q_{\sD}$ invariant; that is,
\beq\label{invariant}
\bU^T \bQ_\sD \bU = \bQ_\sD, \qquad \text{for all } \bU \in \sA,
\eeq
a relation which needs only be checked on the four generators $\bS_i$.
In particular it preserves the level sets
$Q_{\sD}(\mathbf{v}) = k$. It acts on all (real) Descartes 
quadruples $\mathbf{v} = (a,b,c,d)^T$ (viewed as column vectors) via
matrix multiplication, sending $\mathbf{v}$ to $U\mathbf{v}.$ 

The group  $\Aut(Q_{\sD},\zz)$ is 
very large, and can be put in one-to-one correspondence
with a finite index
subgroup of the integer Lorentzian group
$O(3, 1, \zz)$ (using the correspondence \eqn{deslor} in \S2). 
The Apollonian group $\sA$ is
of infinite index in $\Aut(Q_{\sD}, \zz)$, 
and it has infinitely many distinct integral orbits for each integer
$k \equiv 0$ or $3$ $(\bmod~4)$, the values  
for which  integral solutions exist.

Geometrically, the generators $\bS_i$ in $\sA$
correspond to inversions with respect to a circle passing through
three intersection points in
 a Descartes configuration.  This inversion gives a new
Descartes configuration in the same packing, leaving
three of the circles fixed; see \cite[Section 2]{GLMWY11}.
In particular the group 
can be viewed as acting on column vectors $\bv=[a, b, c, d]^T$
giving the ``curvatures'' in a Descartes configuration,
converting  it to the curvatures of another Descartes configuration
in the packing. The generators $\bS_i$
allow us to move
around inside the packing, from one Descartes configuration
to a neighboring configuration. From \eqn{desform} and \eqn{grepn} we see that if 
the curvatures $a,b,c$ of three touching
circles are given, then the curvatures of
the two possible circles which are tangent to these three satisfy
\beqs
d,d' = a + b + c \pm  \sqrt{4(ab+bc+ac)+k}.
\eeqs
Therefore $d + d' = 2(a+b+c)$.  If we begin with circles with curvatures
$a,b,c,d$, the curvature $d'$of the other circle which 
touches the first three is
\beqs
d' = 2(a+b+c) - d. 
\eeqs
The Apollonian group action corresponding
to this sends $(a, b, c, d)$ to $(a, b, c, d')$,
which is the action of $\bS_4$.
In particular, if we start with an integer quadruple,
this procedure will only create other integer quadruples. 
This group is studied in Aharonov \cite{AS97}, 
Graham et al.\ \cite{GLMWY11}, \cite{GLMWY12}
and S\"{o}derberg \cite{So92}.

\begin{defi}~\label{de32}
An {\em integral Descartes ensemble} $\sA[\bv]$ is an orbit of 
the Apollonian group $\sA$ acting on 
an integral quadruple $\bv= (w, x, y, z)^T \in \zz^4$.
That is, 
$\sA[\bv] := \{ \bU \bv:~ \bU \in \sA\}$.
\end{defi}

All elements $\bx \in \sA[\bv_0]$ are integer quadruples satisfying
\beqs
Q_{\sD}(\mathbf{x})= Q_{\sD}(\bv_0)=k.
\eeqs
We are most interested
in  the case $L(\bv) := w+x+y+z \equiv 0~(\bmod~2)$, where $k=4m$. 
For the cases $m= -1, 1$ such ensembles give the set of
Descartes quadruples in an integral spherical (resp.\ hyperbolic)
Apollonian circle packing.  

This section addresses the problem of
classifying the orbits of the
Apollonian group acting on the set of all integral solutions
$Q_{\sD}(\bv) = k$ for a fixed value of $k$. 
%
%
%
%
\subsection{General Reduction Algorithm}

We present  a reduction procedure which, given an integral
quadruple $\bv_0$, finds an element of (locally) ``minimal''  size in 
the orbit $\sA[\bv_0]$. We will show that every orbit contains
at least one and at most finitely many reduced elements.
In most cases the ``minimal''  quadruple is unique and is independent
of the starting point $\bv_0$ of the reduction procedure, and this
is the case for root quadruples defined below. 
However when  $k >0$ there sometimes exist 
exceptional integral Descartes
ensembles containing more than one reduced 
element, and the outcome of the reduction procedure
depends on its starting point in the orbit. 
We show that for each fixed $k$ there are 
in total at most finitely many such exceptional orbits.

The general reduction procedure given below greedily attempts
to reduce the size of the elements in
a Descartes quadruple $\mathbf v := (a,b,c,d)^T$ by applying the
generators $\bS_i$ to  decrease the quantity
\beqs
|\mathbf v| := |a| + |b| + |c| + |d|. 
\eeqs

\noindent {\bf General reduction algorithm.} \\
{\em Input: An integer quadruple $(a,~b,~c,~d)$.

(1) Order the quadruple so that $a \le b \le c \le d.$
Then test in order $1 \le i \le 4$ whether some $\bS_i$
decreases $|\bv| :=|a|+|b|+|c|+|d|.$  For the first $\bS_i$
that  does, 
apply it  to produce a new quadruple, and continue. 

(2) If no $\bS_i$ strictly decreases $|\bv|$, halt.  Declare
the result a {\em reduced quadruple}. 
} \\

This reduction procedure is slightly different from the
reduction algorithm used in part I \cite{GLMWY21}
for the case $k=0$. The part I reduction procedure applied
only to quadruples with $L(\ba)=a+b+c+d >0$, and tried
to decrease the invariant $L(\ba) = a+b+c+d$
at each step,  halting if this could not be done.
Although the general reduction algorithm uses a different
reduction rule, 
one can show in this case that 
it takes the identical sequence
of steps as the reduction algorithm in part I; see 
the remark at the end of \S3.3. 
Thus  the general reduction algorithm can be viewed
as a strict generalization of the algorithm of \cite{GLMWY21}.

We will classify reduced quadruples as either {\em root quadruples} or
{\em exceptional quadruples}, as defined in the
following theorem. 

%
%
%

\begin{theorem} \label{th31}
The general reduction algorithm starting from any nonzero integer
quadruple $\ba= (a, b, c, d)$ always halts in a finite
number of steps. Let the reduced quadruple at termination
be ordered as $\ba= (a, b, c, d)$ with $a \le b \le c \le d.$
Suppose $L(\ba) = a+b+c+d \ge 0$. Then exactly one of the
following holds:

(i) The quadruple $\ba$ satisfies
\beq~\label{nasc1}
a + b+ c \ge d > 0.
\eeq
In this case we call $\ba$ a root quadruple.

(ii) The quadruple $\ba$ satisfies
\beq~\label{nasc2}
a+ b+c \le 0 < d.
\eeq
In this case we call  $\ba$ an exceptional quadruple. 

If $\ba$ is a reduced quadruple such that $L(\ba)< 0$,
then  $\ba^{\ast} := (-d, -c, -b, -a)$ is a reduced quadruple
with $L(\ba^{\ast}) > 0$. We call such $\ba$ 
a root quadruple (resp.  exceptional quadruple)
if and only if 
$\ba^{\ast}$ is a root quadruple (resp.  exceptional quadruple).
\end{theorem}

\paragraph{Proof.}
The integer-valued  invariant
$|\ba| := |a| +|b|+|c| + |d|$
strictly decreases at each step of the reduction algorithm.
It is nonnegative, so the process stops after finitely many
iterations.

Now suppose that $\ba$ is reduced  and
 $L(\ba)= a+b+c+d \ge 0$.   We must
show that (i) or (ii) holds.  Now   $d > 0$ because $\ba$ is
nonzero. The condition that $\ba$ be reduced 
requires $|\bS_4\ba| \ge |\ba|$, and since  
$\bS_4\ba =(a, b, c, d')$
with $d' = 2(a+b+c) -d,$ this condition is  $|d'| \ge |d| = d$. 
Now $d' \ge d$ gives $a+b+c \ge d$ which is (i),
while  $d' \le -d$ gives $a+b+c \le 0$ which is (ii). 
~~~$\Box$ \\

The definition of root quadruple formulated in Theorem~\ref{th31}
is not identical in form to the definition of
root quadruple used in part I \cite{GLMWY21} for 
the case $k=4m=0$. However they are equivalent definitions, since
the reduction algorithms take identical steps, 
as remarked above.

%
%
%
%

\subsection{Root Quadruples}

We first show   that the reduction algorithm is 
well behaved for all orbits $\sA[\bv]$ containing a root quadruple.

%
%
%

\begin{theorem}~\label{th32}
Let $\sA[\bv]$ be a nonzero  integer orbit of
the Apollonian group, and suppose that $\sA[\bv]$
contains a root quadruple $\ba$. Then:

(1) The quadruple $\ba$ is the unique reduced quadruple in 
 $\sA[\bv]$.

(2) All values $L(\bx)$ for $x \in \sA[\bv]$
are nonzero and have the same sign. If this sign
is positive, then the general reduction algorithm starting
from any $\bx \in \sA[\bv]$ applies only $\bS_4$ until the
root quadruple is reached.  
\end{theorem}

\begin{pf}
It suffices to prove the 
result in the case when the  root quadruple has
$L(\ba) = a +b + c+d \ge 0$,
because the case $L(\ba) < 0$ can be treated  using $\ba^{\ast}$. 
Without loss of generality we may
reorder the quadruple  $a \le b \le c \le d$
(since $\sA$ is invariant under conjugation by elements of
the symmetric group on $4$ letters, represented as permutation matrices).
 Now the root quadruple condition \eqn{nasc1} yields  
\begin{equation}~\label{eq31a}
L(\ba) = a+b+c+d >  a+b+c \ge d > 0.
\end{equation}
This implies that
\begin{equation}~\label{eq31b}
a+b \ge d - c \ge 0,
\end{equation}
where the last inequality follows from the ordering.\\

{\bf Claim.} {\em Let $\bb \in  \sA[\bv]$ be
a Descartes quadruple and let $\bP_{\sigma}\bb = (w, x, y, z)$
be its coordinates permuted into increasing order
$w \le x \le y \le z$. Then these coordinates satisfy:
\begin{equation*}~\label{P1}
\mbox{Property (P1):}~~~~~~~~~~ w+x \ge 0.
\end{equation*}
If $\bb$ does not equal the given root quadruple $\ba$ 
then these coordinates  satisfy:
\begin{equation*}~\label{P2}
\mbox{Property (P2):}~~ 0 <  w + x + y <  z.
\end{equation*}
}

To  prove the claim, we observe that every 
quadruple $\bb \ne \ba$ in
$\sA[\bv]$ can be written as
$\bb = \bS_{i_k}\bS_{i_{k-1}} \cdots \bS_{i_1} \ba$ with
each $\bS_{i_j} \in \{ \bS_1, \bS_2, \bS_3, \bS_4\}$, for  some $k \ge 1$,
and with $\bS_{i_{j}} \ne \bS_{i_{j-1}}$ for $2 \le j \le k$,
since all $\bS_i^2= I$. 
We proceed by induction on the number of multiplications
$k$. 

The case $k=0$ is the root quadruple $\ba$,
which satisfies property (P1) by \eqn{eq31b}.
For the base case $k=1$ we must verify
that  properties (P1) and (P2)
hold for
 $\bS_j \ba$ for $1 \le j \le 4$,
where $\ba$ is the given root quadruple. Without loss of
generality we may assume $\ba=(a,b,c,d)$ has the ordering
$a\le b \le c \le d$; a permutation of
$\ba$ merely  permutes the four values   $\bS_j \ba$.
Consider first
$\bS_1 \ba = (a', b,c, d)$. We assert that the
ordering of the elements of $\bS_1 \ba$ is
\begin{equation}~\label{eq31c}
b \le c \le d < a'.
\end{equation}
To establish this we first show 
\begin{equation}~\label{eq31d}
a' := 2(b+c+d) - a =   (b+c+d) + (b+c+d - a) >  b+c+d. 
\end{equation}
This inequality holds since $b - a \geq 0$ from the ordering
while $d > 0$ and $c \geq 0$ from \eqn{eq31a}.
This yields  $ a' > b+ c+d \ge (a+b) +d \ge d,$
which establishes \eqn{eq31c}.
Property (P1) for $\bS_1\ba$  asserts that  $b+c \ge 0$, 
and this holds by \eqn{eq31b}. Property (P2) 
for  $\bS_1\ba$  asserts that
$0< b+c+d < a'$, and 
the left inequality follows from
$b+c +d \geq (a + b) + d >0$, using both \eqn{eq31a} and \eqn{eq31b},
while the right inequality follows from \eqn{eq31d}.
This proves the base case for $\bS_1 \ba$.
The proofs of the base cases $\bS_j\ba$ for $j=2, 3$
are essentially the same as for  $\bS_1\ba$. In these cases the
one curvature that changes becomes strictly larger than $d$.
Finally we treat   $\bS_4 \ba = (a, b, c, d')$.
The ordering that holds here is
\beqs
a \le b \le c \le d'.
\eeqs
Indeed we have
\beq\label{eq:3.11}
d'= 2(a+b+c) - d =  (a+b+c) + (a+b+c-d) \ge a+b+c \ge d,
\eeq
by \eqn{eq31a}. If
 $d=d'$  we stay at  the root
quadruple, so we may assume $d' > d$ in what follows.
Property (P1) for $\bS_4\ba$  asserts $a+b \ge 0$, and this holds
by \eqn{eq31b}. 
Property (P2) for $\bS_4 \ba$ asserts
$0 < a+b+c < d'$.
Now if $ a+b+c=d$ then $d'=d$, so assuming $d' > d$ we
must have $a+b+c> d$.  
Using this fact, we see that the inequalities in 
\eqn{eq:3.11} are both strict, which gives the right inequality 
$a + b + c < d'$.
The left inequality $a+b+c >0$ holds by \eqn{eq31a}.
Thus properties (P1) and (P2) 
hold for $\bS_4 \ba$  as well, or else we stay at the
root quadruple. This verifies (P1) and (P2)
in all cases, and shows that if  $\bb = \bS_j \ba \ne \ba$,
then the largest coordinate in $\bb$ is unique and is the coordinate
in the $j$-th position.

For the induction step, for $k+1$ given $k$, we 
suppose that $\bb = \bS_{i_k}\bS_{i_{k-1}} \cdots \bS_{i_1} \ba$ with 
$\bb \ne \ba$, $k \ge 1$, and we include in
the induction hypothesis the assertion that 
the largest coordinate of $\bb$ is unique and occurs in 
the $i_k$-th coordinate position.   
Let $\bb'= (w, x, y, z)= \bP_{\sigma} \bb$ 
denote its rearrangement in increasing order $w \le x \le y \le z$,
where $\bP_{\sigma}$ is a $4 \times 4$ permutation matrix.
The induction hypothesis asserts that (P1) and (P2) hold for $\bb'$,
together with the strict inequality  $y < z$.
The induction step requires showing  that properties (P1) and (P2)
hold for $\bS_j\bb'$ for $1 \le j \le 3$, corresponding
to a product of length $k+1$, and that the new coordinate formed
 in $\bS_j\bb'$ is its strictly largest coordinate.
The case of multiplication by  $\bS_4$ is excluded 
because it corresponds to 
choosing $\bS_{i_{k+1}}= \bS_{i_k}$. Indeed
$\bS_4$  changes the largest coordinate, which
$\bS_{i_k}$ did at the previous step by the induction hypothesis,
so that  
$\bS_{i_{k+1}}:= \bP_{\sigma}^{-1} \bS_4 \bP_{\sigma} = \bS_{i_k}.$
The first case is   $\bS_1\bb' = (w', x,y, z).$ We assert
that 
\begin{equation}~\label{eq31e}
x \le y \le z < w'.
\end{equation}
Indeed $x,y,z \ge 0$ by property (P1) for $\bb'$ and 
$$
w' = 2(x+y+z) - w = (x+y+z) + (x+y+ z -w) > x+y+z,
$$
since property (P2) for $\bb'$ gives $x+y+z > w+2x+2y \ge w.$
Next property (P1)  for $\bb'$ gives 
$$
x+y+z \ge (w+x)+z\ge z > y \ge 0,
$$
whence  $w' > x+y+z \ge z$ and \eqn{eq31e} follows.
Property (P1) for $\bS_1\bb'$ asserts $x+y \ge 0$, which holds
since $x+y \ge w+x \ge 0$, using Property (P1) for $\bb'$.  
Property  (P2) for $\bS_1\bb'$ asserts 
$0< x+y+z < w'$, both inequalities of which were shown above.
The second case is  $\bS_2\bb' = (w, x', y, z)$. We assert that
\begin{equation}~\label{eq31f}
w \le y \le z < x',
\end{equation}
where we must  show $x' > z$. Now
$x' = 2(w+y+z)-x = (w+y+z) + (w+y+z-x) > w+y+z$
since property (P2) for $\bb'$  gives $w+y+z > 2w+ x+2y \ge x$.
Next $w+y \ge w+x \ge 0$ holds by
property (P1) of $\bb'$, so $w+y+z \ge z >0$
whence $x' > z$, and \eqn{eq31f} follows. 
Property (P1) for $\bS_2\bb'$ now asserts
$w+y \ge 0$, which holds since $w+y+z \ge z$.
Property (P2) for $\bS_2\bb'$ now  asserts $0< w+y+z < x'$,
and both of these inequalities were verified above.
Finally, the  proof for $\bS_3\bb'=(w, x, y', z)$ is similar,
where one shows  $w \le x \le z < y'$.
This completes the induction step, 
and the claim follows. 

To continue the proof,
we observe that for any non-root quadruple $\bb'$
 property (P2) implies that  
$L(\bb) = (w+x+y)  + z > 0$. Thus all quadruples in
the packing have $L(\bb) > 0$. We next show $\bb$
is not reduced.  Indeed property (P2) gives
$\bS_4 \bb = (w, x, y, z')$
with $ 0< z' = 2(w+x+y) - z < z$ so that
$|\bS_4 \bb| < |\bb|$. 
We conclude that  the root quadruple is the
unique reduced quadruple, which proves assertion (1).

To verify assertion (2), we must show that
$|\bS_j\bb| \ge |\bb|$ for $1 \le j \le 3$, while
$|\bS_4 \bb| < |\bb|.$ The last inequality is already done.
Now $\bS_3 \bb = (w, x, y', z)$ with $y' = 2(w+x+z) -y \ge 
2(w+x) + 2(w+x+y) -y \ge y$ has $|\bS_3\bb| \ge |\bb|$.
The argument showing  $|\bS_j \bb| \ge |\bb|$ if $j=1, 2$ 
is similar. This proves (2). 
~~~$\Box$
\end{pf} \\

For $Q_{\sD}(\ba)= k$ with $k$ fixed there are infinitely
many distinct root quadruples.
The next result determines  various of their properties.

%
%
%
%

\begin{theorem} \label{th33}
(1) Any root quadruple $\ba = (a,b,c,d)$ with 
$Q_{\sD}(\ba)=k$ and 
$a \leq b \leq c \leq d$ and $L(\ba) > 0$
satisfies
\beq~\label{scond}
0 \le b \le c \le d.
\eeq
In addition, if $k > 0$, then 
\beq~\label{scond2}
a < 0,
\eeq
and if $k \le 0$ then 
\beq~\label{scond3}
a \le \sqrt{|k|}.
\eeq

(2) For each pair of integers $(k,-n)$
there are only finitely
many root quadruples $\ba$ with 
$Q_{\sD}(\ba) = k$ and $a=-n$, with
the exception of those pairs $(k, -n)=(4l^2, -l)$ with
 $l \ge 0$. In 
each latter case there is an infinite family of root
quadruples $\{ (-l, l, c, c): ~ c \ge \max(l, 1) \}.$
\end{theorem}

\begin{pf}
(1)  The condition \eqn{nasc1} implies
$a+b \ge d- c \ge 0$, so that $b \ge 0$, and the ordering
gives \eqn{scond}. 

We view the Descartes equation $Q_{\sD}(\ba) = k$ 
as a quadratic equation in $a$, and solving it gives
\beqs
a = b+c +d \pm \sqrt{4(bc + bd + cd) + k}
\eeqs
The ordering and \eqn{scond} give 
$a \le b \le b + c + d$, which shows that the minus sign
must be taken in the square root, and
\beq~\label{scond5}
a = b+c +d - \sqrt{4(bc + bd + cd) + k}
\eeq
Using \eqn{nasc1} we have
\beqs
2(bc + bd + cd) = (bc + d(b+c)) + (bc+bd+cd) \\
\ge (b^2 + d(d-a)+ (bc +bd + c^2) \\
 = b^2 + c^2 + d^2 + [bc + d(b-a)],
\eeqs
which gives
\beqs
4(bc + bd + cd) \ge  (b+c+d)^2 + [bc + d(b-a)].
\eeqs
As a consequence \eqn{scond5} yields
\beq~\label{scond6}
a \le (b+c+d) - \sqrt{ (b+c+d)^2 +k}.
\eeq
For $k > 0$ this immediately gives
\beqs 
a \le  b+c +d - \sqrt{(b+c+d)^2 +1} < 0,
\eeqs
which is \eqn{scond2}.

For $k \le 0$ we must have
$w := b+c+d \ge \sqrt{|k|}$ so that the
square root above is real, whence \eqn{scond6} gives
\beqs
a \le \max_{w \ge \sqrt{|k|}}[ w - \sqrt{w^2 +k}] = \sqrt{|k|},
\eeqs
which is \eqn{scond3}. 
Indeed the function $f(w)= w - \sqrt{w^2 +k}$
on this domain  is maximized at its left
endpoint $w= \sqrt{|k|}$, since
\beqs
f'(w) = 1 - \frac{w}{\sqrt{w^2 +k}} < 0 ~~\mbox{for}~~ w > \sqrt{|k|}.
\eeqs

Root quadruples with $a>0$ 
do  occur for some negative $k$.
For example $\bv= (2, 3, 4, 7)$ is a root quadruple
for $k= Q_{\sD}(\bv) = -96.$  

(2) We take the values $Q_{\sD}(\ba) =k$
and $a= -n$ as fixed. To show finiteness, it suffices to bound
$c$ above, in terms of $k$ and $n$.
 The conditions  $a+ b +c \ge d$ and
$0\le b \le c \le d$ yield   $0 \le d \le 3c$,
$0 \le b \le c$ and $-2c \le a \le c$. 

We have  $d=a + b + c -x$ for some $x \ge 0$,
and substituting  this  in \eqn{scond5} gives
\beqs
a= a+ 2b +2c - x - \sqrt{4(b+c)(a+b+c -x) + 4bc +k}.
\eeqs
This simplifies to  
\beqs
0 = 2b + 2c  + x - 
\sqrt{(2b+2c - x)^2 + \left(-x^2 + 4ab +4ac + 4bc +k\right)},
\eeqs
so we must have 
\beq~\label{arf}
x^2 = 4ab + 4ac + 4bc +k.
\eeq
Now $ d \ge c + (a+b -x) \ge c$ gives 
\beqs
a+ b \ge x \ge 0. 
\eeqs
Squaring the left inequality  and rearranging  gives 
$a^2 + 2ab + b^2 \ge x^2  = 4ab+  4ac + 4bc +k$, which yields
\beq~\label{badstuff}
(a-b)^2 \ge 4(a+ b)c +k 
\eeq 
We have $ c \ge b \ge 0$, and holding $a$ and $k$ fixed and letting
$b$ grow, one sees that the right side 
grows at least as fast as $4b^2$ while the left side grows
like $b^2$, so one concludes that 
$b$ is bounded above. One can 
now check that \eqn{badstuff} implies 
\beqs
b \le 3|a| + |k|.
\eeqs
This bound shows that
the left side of \eqn{badstuff} is bounded above, yielding
\beqs
4(a+b) c \le (4|a|)^2 + |k|.
\eeqs
We have $a+b \ge x \ge 0$, and 
whenever $a+b > 0$ then this  inequality bounds $c$ above, 
with $ c \le (|a|)^2 + |k|$,   as was
required. 

There remains the case  $a+b=x= 0$. Letting $a=-l$, we have
$b=l$ and $d= a+b+c - x = c$. Thus we have the quadruple
$\ba= (-l, l, c, c)$ with the requirement $c \ge l \ge 0$.
All of these are root quadruples except for $l=c=0$ which
gives the excluded value $(0,0,0,0)$. Here
$Q(\ba) = 4l^2$, and for each $l \ge 0$ we obtain the
given infinite family of root quadruples. ~~~$\Box$
\end{pf}

\paragraph{Remark.} For $k \ge 0$ the
 upper bound $a \le \sqrt{|k|}$ in Theorem~\ref{th33}(1) 
is not sharp in general. For the 
case $k=-4$ of spherical root quadruples, this upper bound 
gives $a \le 2$, but all spherical root quadruples satisfy 
the stronger bound $a \le 0$. To see this, observe that
$a=2$ can only hold
with  $b+c+d=2$ in \eqn{scond6}, and since $b \le c \le d$ are nonnegative
integers we must have $b \le 0$ so $a \le 0$, a contradiction.
Similarly $a=1$ gives $1 \le b+c+d \le \frac{5}{2}$ in \eqn{scond6},
with  the same contradiction.  The case $k=0$ does occur
with the spherical root quadruple $(0, 1, 1, 2)$ pictured in
Figure~\ref{fig:spherical}.

\paragraph{Remark.} We conclude this subsection by showing 
that the notions of reduction algorithm and root quadruple used
here agree for $k=0$ with those used in part I \cite{GLMWY21}.
The part I reduction procedure applied
only to quadruples with $L(\ba)=a+b+c+d >0$, and tried
to greedily decrease the invariant $L(\ba) = a+b+c+d$
at each step,  halting if this could not be done.
Although the  general reduction algorithm uses a different
reduction rule, in this special case it
 takes the identical series of steps as the reduction 
algorithm in part I, 
so defines the identical notion of root quadruple.
Indeed, suppose  $k=0$ and $L(\ba) > 0$.
We will need a result from 
Theorem~\ref{th34}(3) below, which shows that there are no
exceptional quadruples for $k=0$. so the general reduction
algorithm always halts at a root quadruple.
Comparing  Theorem~\ref{th32}(2) with Lemma 3.1 (iii) of \cite{GLMWY11}
shows all the steps are identical.  Namely, after reordering the
curvatures in increasing order, both algorithms
always  apply $\bS_4$. Finally we check that the halting rules
coincide, i.e., that $|\ba|$ is minimized where $L(\ba)$ is
minimized. For both algorithms the minimal element in
a root quadruple is non-positive, with the other three
elements nonnegative (see Theorem~\ref{th33}(1) and  
\cite[Lemma 3.1(3)]{GLMWY11}). 
Furthermore the first non-positive
element encountered in either algorithm is never changed
subsequently. Since the two algorithms coincide up
to this point, and afterwards the two invariants remain
in the fixed relation   $|\ba| = L(\ba) + 2|a|$,
for fixed  $|a|$, their halting criteria coincide.

%
%
%
%

\subsection{Exceptional Quadruples}

We next consider orbits $\sA[\bv]$ that contain
exceptional quadruples.

%
%
%
%
\begin{theorem}~\label{th34}
Let $\sA[\bv]$ be a nonzero  integer orbit of
the Apollonian group. 

(1) If $\sA[\bv]$ contains
an exceptional quadruple, then it contains elements
$\bv_1, \bv_2$ with $L(\bv_1) > 0$ and $L(\bv_2) < 0$. 

(2) There can be more than one exceptional quadruple in
$\sA[\bv]$.

(3) Exceptional quadruples can exist only for $k \ge 1$.
For each $k \ge 1$ there are only finitely many exceptional
quadruples with $Q_{\sD}(\ba) =k$, all with $H(\ba)^2 \le 2k^2$.  
\end{theorem}

\paragraph{Proof.}
(1) Let $\ba=(a, b,c,d)$ be an exceptional quadruple
in $\sA[\bv]$, with $a \le b \le c \le d$. It suffices
to treat the case $L(\ba) = a + b+c +d \ge 0$, since we can
otherwise apply the argument to $\ba^{\ast}$. 

Now Theorem~\ref{th31} gives $a+b+c \le 0$ and $d > 0$.
Take $\bv_1= \bS_1\ba = (a', b, c, d)^T$ and
$L(\bv_1) = 3(a+ b+c+d) - 4a >0$ if $a < 0$, while
if $a=0$ then $a=b=c=0$ and so $d >0$ and $L(\bv_1)>0$
in this case. Next, take $\bv_2= \bS_4\ba= (a, b,c ,d')$.
Then $L(\bv_2) = 3(a+ b+c) - d <0$ if $d>0$, while
if $d=0$ then $a=b=c=d=0$, a contradiction. 

(2) The quadruple  $\ba= (-2, -3, -5, 12)$ with $|\ba|= 22$ 
and $Q_{\sD}(\ba) = k= 360$ is an exceptional  quadruple. 
Its neighbors are 
$ (-2, -3, -5, -32),  (-2, -3, 9, 12), (-2, 7, -5, 12),$
and $ (6, -3, -5, 12)$. Now $ (-2, -3, 9, 12)$
has a neighbor $\bb= (-2, -3, 9, -4)$ which is also an exceptional quadruple.


(3) Let $k$ be fixed.
Let $\ba= (a, b,c, d)$ be an exceptional quadruple
ordered with 
$a \le b \le c \le d$.
By hypothesis
\beq~\label{eq330}
|\bS_j(\ba)|^2 \ge |\ba|^2 ~~\mbox{for}~~ 1 \le j \le 4.
\eeq
We will show that necessarily $ k \ge 1$ and that
\beq~\label{eq331}
H(\ba)^2 = a^2 + b^2 + c^2 + d^2 \le 2k^2.
\eeq

The minimality property \eqn{eq330} and the bound \eqn{eq331} are
both invariant under reversing all signs, so 
by multiplying by $-1$ if necessary we may suppose 
that $L(\ba) = a+b+c+d \ge 0$
without loss of generality. 
 Since $(a,b,c,d) \ne (0,0,0,0)$
we have $d > 0$. The condition $|\bS_4 \ba|^2 \ge |\ba|^2$
requires that either (i) $a+b+c \ge d$ or  (ii) $a+b+c \le 0$.
By Theorem~\ref{th31}
exceptional quadruples 
correspond to case (ii), and those with $L(\ba) >0$ are  
characterized by the condition
\beqs
d > 0 \ge a +b + c.
\eeqs
This condition necessarily implies $a \leq 0$. We consider 
three exhaustive cases.
In each case, we change notation so that $a,b,c$ are nonnegative;
for example, if $a \leq 0 \leq b \leq c \leq d$, 
then we multiply $a$ by $-1$ so that $\ba = (-a, b, c, d)$.

{\em Case 1. $-a \le -b \le -c \le 0 < d.$}

We have 
\begin{eqnarray*}
k  =  Q_{\sD}(\ba) 
& = & a^2 + b^2 +c^2 +d ^2 -2ab - 2ac -2bc + 2ad +2bd + 2cd \\
 & = & a^2 + b^2 +c^2 +d^2 + 2b(d - a) + 2a(d-c) + 2c(d-b) \\
&\ge & a^2 + b^2 +c^2 +d^2 = H(\ba)^2.
\end{eqnarray*}
Since $\ba \ne (0, 0, 0, 0)$ we infer $k \ge 1$, and 
also \eqn{eq331} holds in this case.

{\em Case 2.  $ -a \le -b \le  0  \le c \le d$.} 

We have
\begin{eqnarray*}
k  =  Q_{\sD}(\ba) & = &  
a^2 + b^2 +c^2 +d ^2 -2ab +2ac +2bc + 2ad +2bd - 2cd \\
&= & (a-b)^2 +(c-d)^2 +2ac +2bc + 2ad +2bd. 
\end{eqnarray*}
Every term on the right is nonnegative, which implies $k \ge 0$.
This gives the additional information 
(since $a \ge b$ and $d \ge c$) that
\beqs
b \le a \le b + \sqrt{k} 
~~\mbox{and}~~ c \le d \le c+ \sqrt{k}.
\eeqs
as well as $ac, bc, ad, bd \le \frac{k}{2}$. 
If either of $a, b$ is nonzero
then each of $c, d$ is at most $\frac{k}{2}$. Similarly if either 
of $c, d$ are nonzero then $a, b$ are at most $\frac{k}{2}$. 
Now suppose  both  $a=0$ and  $b=0$.  Then $k = (c-d)^2$,
and there is a reduction step that takes
$(0, 0 , c, d) \mapsto (0, 0, c, 2c-d)$, whence $|2c-d| \ge d$,
so that either $c = d$, or $c=0$. The case $c=d$ gives a root quadruple,
which is excluded, and in the remaining case  $(0,0,0, d)$ we
have $k = d^2 > 0$, and then  $c, d \le \sqrt{k}$,
and in all these subcases
\beqs
H(\ba)^2 = a^2 + b^2 +c^2 +d^2 \le k^2
\eeqs
holds. The argument if $c=d=0$ is similar. In all these cases
we infer that $k > 0$, and that \eqn{eq331} holds.

{\em Case 3. $ -a \le 0 \le b \le c \le d.$} 

We have
\begin{eqnarray}~\label{case3a}
k = Q_{\sD}(\ba) 
& = & a^2 + b^2 +c^2 +d ^2 + 2ab +2ac -2bc + 2ad -2bd - 2cd 
\nonumber \\
&= &(a+b)^2 + (c-d)^2 + 2 (a-b) (c+d)
\end{eqnarray}
In the exceptional case we must have $-a + b+ c \le 0$,
so that $a \ge b+c$. This gives $a \ge b$, whence all terms on
the right side of \eqn{case3a} are nonnegative, and one is strictly
positive so $k > 0$. Also \eqn{case3a} gives $a, b \le \sqrt{k}$,
and using   $a \ge b+c$ it also gives
\beqs
k \ge (c-d)^2 + 2c(c+d) \ge 3c^2 + d^2.
\eeqs
It follows that 
\beqs 
H(\ba)^2 = a^2 + b^2 +c^2 +d^2 \le 2k \le 2k^2,
\eeqs
and \eqn{eq331} holds. 
~~~$\Box$ \\

On comparing Theorem~\ref{th33} and Theorem~\ref{th34}
we see that: {\em  An integer Apollonian group orbit $\sA[\bv]$
 contains an exceptional
quadruple if and only if it contains quadruples with
$L(\bx_1)>0$ and $L(\bx_2) < 0$.}

We conclude this section by determining  the
 exceptional quadruples
for spherical and hyperbolic Apollonian packings,
i.e., for $k=  \pm 4$. There are no exceptional
quadruples for 
spherical packings ($k=-4$) by Theorem~\ref{th34} (3),
so we need only treat the hyperbolic case. 

%
%
%
\begin{theorem}~\label{th35}
For $k=4$ the exceptional orbits $\sA[\bv]$
are exactly those 
containing an element that has a zero coordinate.
There are exactly
two such orbits, one containing  
two exceptional quadruples $(0,0,0, \pm2)$,
and the other containing the exceptional quadruple
$(-1, 0,0, 1)$.
\end{theorem}

\begin{figure}
\begin{center}
\includegraphics[height=3.5in]{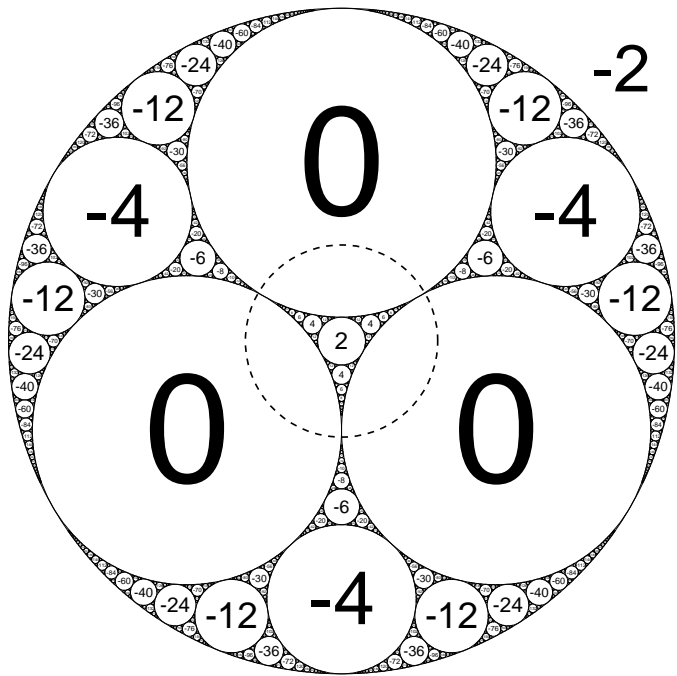}
\end{center}
\caption{The exceptional hyperbolic Apollonian circle packing (0,0,0,2)} 
\label{fig:exception1}
\end{figure}

\begin{figure}
\begin{center}
\includegraphics[height=3.5in]{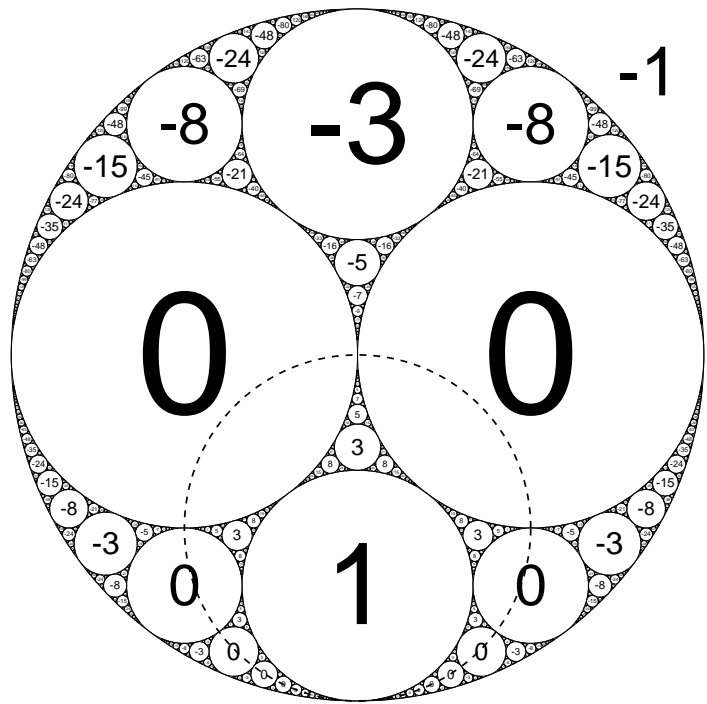}
\end{center}
\caption{The exceptional hyperbolic Apollonian circle packing (-1,0,0,1)}
\label{fig:exception2}
\end{figure}

\begin{pf}
By Theorem~\ref{th32}(3), there are finitely many exceptional
quadruples for $k = 4$.  From \eqn{eq331}, we see that if $(a,b,c,d)$
is an exceptional quadruple, then $a^2 + b^2 + c^2 + d^2 \leq 32$.  A
computer search shows that there are exactly 7 
quadruples with $a \leq b \leq c \leq d$ satisfying this bound along with
$a + b + c + d \geq 0$.
They are
\beqs
\{(-1,0,0,1), (0,0,0,2),
(-1,1,1,1), (-1,1,2,2),
(-1,1,3,3), 
 (0,0,1,3), (0,0,2,4) \}.
\eeqs
It is easy to check that the first five are all reduced while the last two 
are not.  Furthermore, the  third, fourth and 
fifth quadruples are  root quadruples, so may be discarded. 
Thus only  $(-1,0,0,1), (0,0,0,2)$ remain, and these are exceptional.
The orbit of $(0,0,0,2)$ also contains the reduced quadruple $(-2,0,0,0)$; 
see Figure~\ref{fig:exception1}.
The packing with exceptional quadruple $(-1,0,0,1)$ is shown in
Figure~\ref{fig:exception2}.

If fact, if $\sA[\bv]$ does contain a circle of curvature 0, then
it necessarily contains an exceptional quadruple.
Setting $a=0$, the orbit contains a 
integer quadruple $(0,b,c,d)$ with
\beq \label{eq:F}
F(b,c,d) := b^2 + c^2 + d^2 - 2(bc + bd + cd) = 4.
\eeq
Without loss of generality, 
assume that the elements are ordered $|b| \leq |c| \leq |d|$ and
that $b \leq 0$.
We claim that necessarily $b = 0$.

{\em Case 1.} $b \leq 0$, $c > 0$, $d > 0$.\\
Now 
\beqs
f(b,c,d) &=& b^2 + c^2 + d^2 + 2|bc| + 2|bd| - 2|cd| \\
&=& b^2 + (c-d)^2 + 2|bc| + 2|bd|
\eeqs
If $|b| \geq 1$, then $|c|, |d| \geq 1$ also and
$F(b,c,d) \geq 1^2 + 0^2 + 2 + 2 = 5$
which gives us a contradiction.  So we must have $b = 0$.

{\em Case 2.} $b \leq 0$, $c \leq 0$, $d > 0$.\\
Here $|b| \leq |c| \leq |d|$ gives $|bd| > |bc|$, 
so
\beqs
F(b,c,d) &=& b^2 + c^2 + d^2 + 2|bd| + 2|cd| - 2|bc| \\
&\geq& b^2 + c^2 + d^2 + 2|cd|
\eeqs
If $|b| \geq 1$, then so are $|c|$ and $|d|$, so $F(b,c,d) \geq 5$, a
contradiction.  Thus $b=0$.

{\em Case 3.}  $b, c, d \leq 0$.\\
First, apply the reduction algorithm to $(0,b,c,d)$.  
Notice that the first coordinate
always stays $0$ because $\bS_1$ can never make it any smaller.
The algorithm halts when $\bS_4$ can no longer reduce 
$|b| + |c| + |d|$, that is,
when  $|d'| \geq |d|$.  Combining this with $d + d' = 2(b + c)$, we get 
$|d| > |b + c|.$  Therefore
\beqs
b^2 + c^2 + d^2 \leq b^2 + c^2 + (b+c)^2 = 2b^2 + 2bc + 2c^2.
\eeqs
Since $bc, bd, cd \geq 0$, we have
$bc + bd + cd \geq |b|^2 + |bc| + |c|^2.$
Thus,
\beqs
F(b,c,d) = b^2 + c^2 + d^2 - 2(bc + bd + bd) \leq 0
\eeqs
and hence cannot satisfy \eqn{eq:F}, a contradiction.  So there are no
solutions in this case.

The claim $b=0$ is proved, 
so any quadruple with one element zero is contained in 
the same packing as some $(0,0,c,d)$.
Now solve
$2(c^2 + d^2) =  (c + d)^2 + 4$
and we see that we must have $(0,0,c,d) = (0,0,n,n+2)$ for some
$n \in \zz$.  Reduction by the Apollonian group changes 
$(0,0,n,n+2) \to (0,0,n,n-2)$, and therefore there are only
two packings containing a curvature zero circle, those are the
packing containing $(0,0,0,2)$  and the packing 
containing $(-1,0,0,1)$.~~~$\Box$
\end{pf} \\

To understand the geometric
nature of the exceptional hyperbolic 
packings, we  identify the real hyperbolic plane $\mathbb H$
with the interior of
the unit disk, with ideal
boundary the circle  of radius $1$ centered at the origin.
Recall that the curvature of a hyperbolic circle is given by 
$\coth r$, where $r$ is the hyperbolic radius of the circle.  
For circles contained in the real hyperbolic plane, $r$ is positive and
\beqs \label{3.6}
\coth r = \frac {e^r + e^{-r} } {e^r - e^{-r}} \ge 1
\eeqs 
with equality if and only if $r = \infty$. 
However we  will  allow 
hyperbolic Apollonian circle
packings which lie in the entire plane. 
We have a second copy of
the hyperbolic plane which is the exterior of the closed unit disk,
and we use the convention that (positively oriented)
circles entirely contained in
the exterior region have the sign of their curvature reversed,
 in which case $\coth r$  ranges from
$-1$ to $-\infty$. Euclidean circles  that intersect the 
ideal boundary (unit circle) are treated  as hyperbolic circles
of pure imaginary radius, and these cover the remaining range
$-1 \le \coth r \le 1$. A hyperbolic circle of
curvature $0$ is a hyperbolic
geodesic, which corresponds to a Euclidean circle which intersects
the ideal boundary (unit circle)  at right angles, and a hyperbolic circle of
curvature  $1$ is a horocycle, i.e., a circle tangent to the
ideal boundary. 

The two exceptional integral hyperbolic packings
are pictured in Figures~\ref{fig:exception1} and \ref{fig:exception2}, 
respectively.
The pictures are not drawn to the same scale; in each one 
the dotted circle represents 
the ideal boundary of hyperbolic space, and is not
a circle belonging to the packing. The hyperbolic packing
$(0,0,0,2)$ contains exactly three circles assigned
hyperbolic curvature zero, while the 
hyperbolic packing $(-1, 0,0, 1)$
contains infinitely many circles assigned hyperbolic curvature zero.

%
%
%
%
%

\section{Distribution of Integer Root Quadruples}
\label{root}
\setcounter{equation}{0}

In this section we restrict to the case  $k=4m$.
We enumerate root quadruples $\bv= (w,x,y,z)$, and without
essential loss of generality we only consider such quadruples    
with $L(\bv)=w+x+y+z > 0$.  The
root quadruples are classified according to the size of their
minimal element, say $w=-n$, where we suppose $n >0$.
In the spherical and hyperbolic cases $k = \pm 4$ 
this corresponds geometrically 
to characterizing the distinct integral Apollonian packings
that have an outer enclosing circle of signed curvature $-n$.
 
\begin{defi} 
Let $N_{\text{root}}(k; -n)$ count the number of integer 
root quadruples $\bv=(w,x,y,z)$ to
$Q_{\sD}(\bv) = k$ with $L(\bv)=w+x+y+z > 0 $ and  $w = -n$.
\end{defi}

In \cite[Sec. 4]{GLMWY21} 
the number $N_{\text{root}}(0, -n)$ of 
Euclidean root quadruples was interpreted
as a class number for
binary quadratic forms under a $GL(2, \zz)$ action.
This interpretation was discovered
starting from an exact formula for the
number of (Euclidean) root quadruples having
smallest element $-n$, originally  conjectured
by C. L. Mallows and S. Northshield,
and independently proved by Northshield \cite{North}.
Below we show that a similar class number interpretation
holds generally for the case  $Q_{\sD}(\ba)= 4m$.

We write a binary 
quadratic form of even discriminant $\Delta = 4B^2 - 4AC$ as
\begin{equation}~\label{eq41}
[A,B,C] := AT^2 + 2BTU + CU^2.
\end{equation}
Here we follow the classical notation of Mathews~\cite{Ma1892},
to maintain consistency with the notation of \cite{GLMWY21}.
(The convention of  
Buell~\cite{Bu89} writes this form as $[A, 2B, C]$,
and allows odd middle coefficients.) 
A form is {\em definite} if $\Delta < 0$, and it is then 
{\em positive} if it only represents positive values (and zero).
Positive definite forms are those definite forms
with both $A$ and $C$ positive, which
follows from $-4AC \le 4(B^2 - AC) < 0.$

The standard action on binary forms is an $SL(2,\zz)$-action.
A {\em reduced form} is one that satisfies
\begin{equation}~\label{red1}
2|B| \le A \le C.
\end{equation}
Every positive definite form is equivalent (under the
$SL(2,\zz)$-action) to at least one  {\em reduced form}.
All reduced forms are $SL(2,\zz)$-inequivalent except for
$[A,B,A] \equiv [A,-B,A]$ and $[A,A,C] \equiv [A,-A, C]$.  

There is also a $GL(2,\zz)$-action on the space of binary
quadratic forms.  It defines an 
action of $\begin{pmatrix}1&0 \\ 0 & -1\end{pmatrix}$ on 
binary forms, sending
$[A,B,C]$ to $[A,-B,C]$.  We define the
 $GL(2, \zz)$-reduced forms to be  
the subset of 
$SL(2,\zz)$-reduced forms satisfying
\begin{equation}~\label{red2}
0  \leq 2B \leq A \leq C.
\end{equation}
Then every positive definite form is $GL(2, \zz)$-equivalent 
to a unique  $GL(2, \zz)$-reduced form.

A form is {\em primitive} if $\operatorname{gcd} (A, 2B, C) =1$ and is
{\em imprimitive} otherwise.
Define $\widetilde{h^{\pm}}(\Delta)$ to be the number of 
(primitive or imprimitive) $GL(2,\zz)$-reduced
forms of discriminant $\Delta$.
Let $\widetilde{h}(\Delta)$ be the number
of (primitive or imprimitive) $SL(2,\zz)$-equivalence classes of forms
of discriminant  $\Delta$.
Since each $GL(2,\zz)$-reduced form $[A,B,C]$ is equivalent to 
the $SL(2,\zz)$-reduced form $[A,-B,C]$, 
there holds
\begin{equation}~\label{classno}
\widetilde{h^{\pm}}(\Delta) = 
   \frac 1 2 (\widetilde{h}(\Delta) + \widetilde{a}(\Delta)),
\end{equation}
in which  $\widetilde{a}(\Delta)$ counts the number of 
(primitive or imprimitive) ambiguous classes. Here an {\em ambiguous class}
is an $SL(2, \zz)$-equivalence class containing 
an {\em ambiguous form}, which 
is a (primitive or imprimitive)
reduced form for which $2B$ divides $A$, i.e., reduced forms of shape
$[A, 0, C]$ or $[2B, B, C]$, cf.  Mathews \cite[Art. 159]{Ma1892}.

The usual class number $h(\Delta)$ is defined as the number of
$SL(2,\zz)$-equivalence classes of primitive forms.
Similarly  we let
$a(\Delta)$ count the number of 
primitive ambiguous classes.
If $[A,B,C]$ is a primitive form of discriminant $\Delta$, then
for $\ell > 1$,  $[\ell A, \ell B, \ell C]$ is an (imprimitive) form of 
discriminant $\ell^2 \Delta$.  We have
\beq~\label{imprimitive1}
\widetilde{h}(\Delta) = 
\sum_{\ell^2 \mid \Delta} h\left( \frac {\Delta} {\ell^2} \right ),
\eeq
where we use the convention that $h(n) = 0$ if $n$ is not a 
valid discriminant.
In the same vein we have
\beq~\label{imprimitive2}
\widetilde{a}(\Delta) = 
\sum_{\ell^2 \mid \Delta} a \left( \frac {\Delta} {\ell^2} \right ).
\eeq
For $\Delta = -4m< 0$  one has 
$a(\Delta) = 2^{t+1}$, where $t$ is the number of distinct prime
divisors of $m$; see Mathews ~\cite[Art. 156]{Ma1892}.

The following result gives a  class number criterion counting
root quadruples for an arbitrary Descartes equation
$Q_{\sD} = 4m$. This result applies even when the associated
discriminant $\Delta$ is non-negative, 
provided the definition of reduced form is
given by \eqn{red1} (resp.\ \eqn{red2}), though this
no longer corresponds nicely to $SL(2, \zz)$-equivalence 
(resp.\ $GL(2, \zz)$-equivalence).

%
%
%
%

\begin{theorem}~\label{th41}
Let $n > 0$ be fixed.

(1) The  
ordered root quadruples $\bv= (w, x, y, z)$ 
with $Q_{\sD}(\bv) = 4m$, for fixed $m$, 
having smallest element $w=-n$ 
are in one-to-one correspondence with reduced 
integral binary quadratic forms
$[A, B, C]= AX^2 + 2BXY + CY^2$ (primitive or imprimitive)
having  discriminant $\Delta = -4n^2 + 4m$ 
and  a non-negative middle coefficient.
If $\bv= (w, x, y, z)$ with $w= -n \le  x \le y \le z$
then the correspondence is given by  
\beq~\label{qformula}
[A,B,C] = \left [-n + x,  \frac{1}{2}(-n + x + y - z), -n + y \right ].
\eeq

(2) When $\Delta = -4n^2 + 4m <0$, the number of 
root quadruples to $Q_{\sD} = 4m$ 
with least element $-n$ 
is given by the $GL(2, \zz)$-class number 
\begin{equation}
N_{\text{root}}(4m; -n) = \widetilde{h^\pm}(-4(n^2 -m)).
\end{equation}
\end{theorem}

\paragraph{Proof.}
 Recall from Theorem~\ref{th31} that an 
ordered integer quadruple $\bv= (w, x,y,z)$
with $w \le x \le y \le z$ and $L(\bv) >0$ 
is a root quadruple if and only if
\beq~\label{rcond3}
w+x+y \ge z \ge 0,
\eeq
and Theorem~\ref{th33} shows
that  root quadruples always satisfy  $0 \le x \le y \le z$.
We now write $w= -n$ to suggest that we are mainly concerned
with negative integers. However, for negative $m$ small
positive values of $w$  are permitted, which necessarily
satisfy $w \le \sqrt{|m|},$ by  Theorem~\ref{th33}(1).

The integer solutions $Q_\sD(w,x,y,z) = 4m$ with $w = -n $ 
are in one-to-one correspondence with integer representations of 
$n^2 -m$ by
the ternary quadratic form
\beqs
Q_T(X,Y,Z) := XY + XZ + YZ,
\eeqs
where
\beqs
(X,Y,Z) := \left (\frac{1}{2}(w+x+y-z), \frac{1}{2} (w+x-y+z), 
\frac 1 2 (w-x+y+z) \right).
\eeqs
The congruence condition $w + x + y + z \equiv 0 \pmod 2$ 
for integer quadruples
implies that $X,Y,Z \in \zz$.
This map is a bijection: given an integer solution to 
$Q_T(X,Y,Z) = n^2 -m$, we write
\beqs
(w,x,y,z) :=(-n, n+X + Y, n + X + Z, n + Y + Z)
\eeqs
and an algebraic calculation shows  that $Q_\sD(w,x,y,z) = 4m$.

The root quadruple conditions \eqn{rcond3} above translate 
to the inequalities
\beq~\label{eq:rootxyz}
0 \leq X \leq Y \leq Z.
\eeq

Next, for any integer $M$, the integer solutions $Q_T(X,Y,Z) = M$ are
in bijection with integer representations of $-4M$ by the determinant
ternary quadratic form
\beqs
Q_\Delta (A,\tilde{B},C) := \tilde{B}^2 - 4AC
\eeqs
A solution $(X,Y,Z)$ gives a solution $Q_\Delta(A,\tilde{B},C) = -4M$ under
\beqs
(A,\tilde{B},C) := (X + Y, 2X, X + Z),
\eeqs
and the inverse map is
\beqs
(X,Y,Z) := (\tilde{B}/2, A - \tilde{B}, C - \tilde{B}).
\eeqs
Note that since $-4M = \tilde{B}^2 - 4AC$ is even, 
then $\tilde{B}$ must also be even, so
$B= \tilde{B}/2$ must be an integer.

One sees easily that the inequalities \eqn{eq:rootxyz} are 
equivalent to the inequalities
\beq~ \label{eq:rootabc}
0 \leq 2B \leq A \leq C.
\eeq
Finally, we recognize that the conditions \eqn{eq:rootabc} give a
complete set of equivalence classes of integral binary quadratic forms
of fixed discriminant $-4M$ under the action of $GL(2,\zz)$.

Therefore, taking $M = n^2 -m$, these two steps associate to any 
(ordered) integer quadruple $(-n,x,y,z)$ the binary quadratic form 
\beqs
[A,B,C] = \left [-n + x, \frac{1}{2}(-n + x + y -z), -n + y \right]
\eeqs
of discriminant $D := -4(n^2 -m)$.  In order for the form
to be positive definite, we must have
$n^2> m $. This imposes no constraint on $n$ if $m <0$,
while if $m \ge 0$ we must have $n > \sqrt{|m|}.$
 This form is $GL(2,\zz)$-reduced if and only if $(-n,x,y,z)$ 
is a root quadruple.

Conversely, given  a form of discriminant $D$, we construct a Descartes
quadruple
\beqs
(w,x,y,z) := (-n, n + A, n + C,  n + A + C - 2B)
\eeqs
which is a root quadruple if and only if $[A,B,C]$ is $GL(2,\zz)$-reduced.
 ~~~$\Box$ \\

Theorem~\ref{th41} counts the number of root 
quadruples $N_{\text{root}}(4m; -n)$ with $L(\bv) > 0$ 
having a fixed smallest element $w= -n$ 
in terms of a class number, provided $n^2 > m$. 
This condition is always satisfied  when $m < 0$. For
the case $m=0$ there is one excluded value
$n=0$ corresponding to multiples of the root quadruple
$(0,0,1,1)$, which falls under  Theorem~\ref{th33}(2).
For positive  $m$ and for the remaining
values $-\sqrt{|m|} \le w < 0$ allowed by Theorem~\ref{th33},
the combinatorial correspondence with certain  binary quadratic
forms of discriminant $\Delta = -4(n^2 - m) \ge 0$ still
applies, but the notion of ``reduced form'' 
in \eqn{red1} is not
that that normally used for indefinite binary quadratic forms,
and $N_{\text{root}}(4m, -n)$ is not interpretable as
a class number. 
For hyperbolic Descartes quadruples ($ m=1$), there is one  
allowable value $n=-1$ not interpretable as a class number. 
Indeed there are infinitely many root quadruples when
$n=-1$ by Theorem~\ref{th33}(2), since $(k,n) = (4m^2, -m)$ with $m=1$.

Table~\ref{tab:root}
presents small values of 
$N^S_{\text{root}}(-n) := N_{\text{root}}(-4, -n)$ 
and $N^H_{\text{root}}(-n) := N_{\text{root}}(4, -n)$.
All  hyperbolic root quadruples have $-n \le -1$, according to
Theorem~\ref{th33}(1), and 
all spherical root quadruples have $n \le 0$, as shown  in
the proof of Theorem~\ref{th33}(1).  The value $n=0$ is not shown
in the table; we have $N^{S}_{\text{root}}(0) =1$, given by the
spherical root quadruple $(0, 1,1,2)$.

%
%
%
%

\begin{table}
\begin{center}
  \begin{tabular}{|c|c|c||c|c|c||c|c|c|}
    \hline 
    $n$ & $N^S(-n)$ & $N^H(-n)$  & $n$ & $N^S(-n)$ & $N^H(-n)$ & 
       $n$ & $N^S(-n)$ & $N^H(-n)$\\\hline\hline
    
    1 & 1 & $+\infty$ & 11 & 6 & 12 & 21 & 6 & 24 \\\hline
    2 & 2 & 2 & 12 & 6 & 12 & 22 & 12 & 12 \\\hline
    3 & 2 & 3 & 13 & 8 & 12 & 23 & 16 & 22 \\\hline
    4 & 3 & 4 & 14 & 6 & 12 & 24 & 5 & 24 \\\hline
    5 & 4 & 6 & 15 & 5 & 19 & 25 & 19 & 24 \\\hline
    6 & 2 & 6 & 16 & 9 & 16 & 26 & 16 & 25 \\\hline
    7 & 5 & 8 & 17 & 12 & 18 & 27 & 8 & 24 \\\hline
    8 & 6 & 8 & 18 & 10 & 10 & 28 & 10 & 28 \\\hline
    9 & 3 & 11 & 19 & 10 & 24 & 29 & 14 & 24 \\\hline
    10 & 8 & 9 & 20 & 11 & 20 & 30 & 14 & 30\\ \hline
  \end{tabular}
\end{center}
\caption{$N_{\text{root}}(\pm4, -n)$ for small $n$.}
\label{tab:root}
\end{table}

The next result  uses bounds for class numbers to
derive an  upper bound for the number of 
root quadruples.

%
%

\begin{theorem}~\label{th43}
For all integers $m$ the  number $N_{\text{root}}(4m; -n)$ of 
ordered integer 
root quadruples $(a,b,c,d)$ with $Q_{\sD}(\ba) = 4m$ and 
$a = -n$ satisfies, uniformly  for $n > \sqrt{|m|}$, 
\beqs
N_{\text{root}}(4m; -n) = O\left( n (\log n)( \log\log n)^2 \right).
\eeqs
The implied $O$-constant is independent of $m$.
\end{theorem}

\begin{pf}
Since $\widetilde{h}^{\pm}(\Delta) \le \widetilde{h}(\Delta)$ we have 
\beqs
N_{\text{root}}(4m; -n) \le \widetilde{h}(-4(n^2 - m))
=\sum_{\ell^2~|~-4(n^2 - m)} h \left( \frac{-4(n^2 - m)}{\ell^2}\right ).
\eeqs
in which  $h(D)$ denotes the usual class number for 
primitive forms under 
$SL(2, \zz)$-equivalence.  
Dirichlet's class number formula gives,
for primitive  negative discriminants $D$, that 
\beqs
h(D) = \frac{w(D) \sqrt{|D|}}{2 \pi} L(1, \chi),
\eeqs
where $\chi$ denotes the (primitive or imprimitive)
character attached to the quadratic
field $\qq(\sqrt{D})$, and $w(D)$ denotes the number of roots
of unity in this field, which is at most six. This formula is
valid for class numbers of non-maximal orders as well,
where the $L$-function now uses an imprimitive character.
For primitive characters we  have the bound
\beqs
L(1, \chi) = O (\log |D|);
\eeqs
 see Davenport~\cite[Chap. 14]{Da80},
and Ramar\'{e} \cite{Ram01} for precise estimates. 
For imprimitive
characters, we must correct by a finite  Euler factor
$\prod_{p \in S} ( 1 - \chi(p)\frac{1}{p})$ involving some subset $S$
of  primes
dividing the discriminant. As shown in \cite[Theorem 4.4]{GLMWY21},
this is bounded above by
$\prod_{p \in S} ( 1 +\frac{1}{p})$, and 
can lead to an extra  factor of size 
at most $O(\log\log |D|)$. Writing
$D = D_0 S^2$ with $D_0$ a fundamental discriminant, we have
\begin{eqnarray}~\label{bbound}
\widetilde{h}(D) &=&O\left(  \sqrt{|D_0|} \log |D_0| 
\sum_{d~|~S}  d \prod_{ p~|~d}\left(1 + \frac{1}{p}\right)
\right) \nonumber \\
& = & O \left(  \sqrt{|D_0|S^2} \log |D_0| (\log\log S)^2 \right). 
\end{eqnarray}
This gives
\beqs 
 \widetilde{h}(D)=    O \left(  \sqrt{|D|} \log |D| (\log\log |D|)^2 \right),
\eeqs
valid uniformly for all discriminants $D \le -4.$ 
Here we used the estimate $\sum_{d~|~S} d  \le S \log\log S$.
The result follows on taking $D= -4(n^2-m)$. The $O$-constant
can be taken independent of $m$ since $n^2 - m \le 2n^2$ under
the stated hypotheses.
~~~$\Box$
\end{pf}

\paragraph{Remarks.} 
(1) The upper bound  of Theorem~\ref{th43}
can be strengthened  when $m=0$.
In that case $D=-4n^2$, and one
obtains the stronger upper bound
$O\left( n (\log\log n)^2\right)$
using \eqn{bbound} with  $D_0=-4$ and
$S= n^2$. This improved upper bound also  follows 
from  \cite[Theorem 4.2]{GLMWY21}. 
The latter result counts primitive root quadruples only
and gets a bound $O (n (\log\log n))$. 
Theorem~\ref{th43} counts primitive and imprimitive
root quadruples (with primitivity defined by the auxiliary
binary quadratic form \eqn{qformula}). Summing over
the imprimitive forms produces an extra factor
of $\log\log n$.

(2) In \cite{GLMWY21} a lower bound
$\Omega \left(\frac{n}{\log\log n}\right)$
was obtained for the number of primitive
Euclidean root quadruples, which implies
\beqs
N_{\text{root}}(0, -n) = \Omega\left( \frac{n}{\log\log n} \right). 
\eeqs
We are unable to obtain  a lower bound for 
of equivalent strength for general $k=4m$, 
because it appears to  require 
obtaining good lower bounds for class numbers of 
general imaginary quadratic fields, which is
a difficult problem.  
Using the
Brauer-Siegel theorem (see \cite[pg.\ 328]{La}) one  can obtain 
for fixed $m$ and any $\epsilon > 0$ a 
lower bound 
\beqs
N_{\text{root}}(4m, -n) \ge n^{1 - \epsilon}
\eeqs
valid for all $n \geq n_0(m,\epsilon)$, in which the constant
$n_0(m,\epsilon)$ is not effectively computable. \\

We end this section by raising some questions 
concerning the behavior 
of the  total number of root quadruples with $Q_{\sD}(\bv) = k$ 
having minimal element below a given bound. 
We define the summatory function
\beqs
S(k, T) := \sum_{n = \lfloor \sqrt{|k|}\rfloor + 1}^{T}
 N_{\text{root}}(k, -n),
\eeqs
Theorem~\ref{th43} gives for nonzero $k=4m$ that 
\beq~\label{sumsum}
S(4m, T) = \sum_{n= \lfloor \sqrt{|k|}\rfloor +1}^T 
\widetilde{h}^{\pm}(4(n^2 +m)).
\eeq
One expects these sums  of $GL(2, \zz)$-class numbers
over a quadratic sequence  to grow roughly like $T^{3/2}$.
Can the leading term in  their asymptotics be determined?
The problem can be reduced to considering similar sums  of
$SL(2, \zz)$-class numbers using \eqn{classno} 
to obtain
\beqs
S(4m, T) = \frac{1}{2} \sum_{n =  \lfloor \sqrt{|k|}\rfloor+1   }^T 
\widetilde{h}(-4(n^2 +m)) + 
O \left( T \cdot 2^{\frac{\log T}{\log\log T}} \right),
\eeqs
where we made use of an upper  bound for the number of ambiguous
forms, as follows. If $t$ denotes the number of distinct
prime divisors of $\Delta <0$ then \eqn{imprimitive2} yields
\beqs
\tilde{a}(\Delta)  \le d(\Delta) 2^{t+1}
 \le 2^{O(\frac{\log \Delta}{\log\log \Delta})}.
\eeqs

One can also ask whether $\frac{S(4, T)} {S(-4, T)}$ tends to
a limit as $T \to \infty$, and if so, compare the limiting
value to that in \eqn{ratiobound}. One expects that there
are asymptotically fewer spherical root quadruples 
than hyperbolic root quadruples. 
However, it is not true that 
$N_{\text{root}}(-4, n) \le N_{\text{root}}(4, n)$ always.
For example, $N_{\text{root}}(-4,32) = 23 > 20 = N_{\text{root}}(4,32)$.

%
%
%
%
%

\section{Integers Represented by a  Packing: Congruence Conditions}
\label{congruences}
\setcounter{equation}{0}

In this section we restrict to the cases $k= \pm 4$ corresponding to 
spherical and hyperbolic Apollonian packings. 
We show there are congruence restrictions on the allowed ``curvatures''
modulo 12 of circles in integral spherical  and hyperbolic Apollonian 
circle packings. 

\begin{theorem}~\label{th51}
In any integral spherical  Apollonian circle packing, 
the ``curvatures'' of the circles in the packing omit exactly
three congruence classes modulo $12$.
\end{theorem} 

\begin{pf}
It suffices to classify the unordered Descartes quadruples in
all integral packings. We claim that
these fall into one of
two eighteen-element orbits mod 12, 
which are $\sO$ and $\sO + 6 \pmod {12}$, where 
\beqs 
\sO =
\{
(2,1,1,0), (5,2,1,0), (5,4,2,1), (6,5,4,1), (6,5,5,4), (8,2,1,1), \\
(8,6,1,1), (8,6,5,1), (9,4,2,1), (9,5,4,2), (9,8,2,1), (10,5,1,0),\\
(10,5,5,0),(10,5,5,4), (10,8,5,1), (10,9,5,4), (10,9,8,1), (10,9,8,5)
\},
\eeqs
and the notation $\sO + 6 \pmod {12}$ means to add 6 to each 
coordinate of each 
element of $\sO$ and reduce modulo 12.

The theorem follows from the claim, since the curvatures in $\sO$
omit the classes $3, 7$ and $11~ (\bmod~12)$, while $\sO + 6$ omits
$1, 5$ and $9 ~ (\bmod~12)$.

To check the claim,
it  is easy to verify that there are $212$ possible solutions 
(without respect to order) of the
spherical Descartes equation \eqn{sphdescartes} modulo 12.  
These lie in 19 orbits under the action 
of the Apollonian group.  However, not all of the solutions 
modulo 12 come from integer solutions.

If $\ba = (w,x,y,z)$ is the reduction 
modulo 12 of a solution to the spherical Descartes equation, 
then there must exist $a,b,c,d$ such that
\begin{multline}\label{sphdes}
2((12a + w)^2 + (12b + x)^2 + (12c + y)^2 + (12d + z)^2) =\\
         (12(a + b + c + d) + (w + x + y + z))^2 - 4
\end{multline}
If we reduce this equation modulo $24$, we get
\beqs
2(w^2 + x^2 + y^2 + z^2) \equiv (w + x + y + z)^2 - 4 \pmod{24},
\eeqs
and we see that $\ba$ must also be a Descartes quadruple mod $24$.
This eliminates 108 of the 212 possibilities and 11 of the 18 orbits.
Furthermore, if we reduce \eqn{sphdes} modulo 48, we get 
the condition that
\begin{multline*}
2(w^2 + x^2 + y^2 + z^2) + 24(a + b + c + d) (w + x + y + z)\\
 \equiv (w + x + y + z)^2 - 4 \pmod{48}.
\end{multline*}
But since $w + x + y + z \equiv 0 \pmod 2$, this is
\beq
2(w^2 + x^2 + y^2 + z^2) \equiv (w + x + y + z)^2 - 4 \pmod{48}.
\eeq
Thus $(w,x,y,z)$ must also be a quadruple modulo 48, 
ruling out 68 of the 104 remaining quadruples.
Finally, a short computer search turns up integer spherical quadruples
in all 36 classes. ~~$\Box$
\end{pf}

\begin{theorem}~\label{th52}
In any integral  hyperbolic Apollonian circle packing, 
the ``curvatures'' of the circles in the packing omit
at least three congruence classes modulo $12$.
\end{theorem} 

\begin{pf}
We classify the (unordered) Descartes quadruples
in an integral hyperbolic packing, showing that they fall in one of 81 
equivalence classes modulo $12$.
We claim that the integer  quadruples lie in 15 orbits under the action of 
the Apollonian group, given by
$\sO_3, \sO_3 + 3, \sO_3 - 3,~~ -\sO_3, -\sO_3 +3, -\sO_3 -3, ~~
\sO_4, \sO_4 + 6, ~~ -\sO_4, -\sO_4 + 6, ~~ \sO_7, \sO_7 + 3, \sO_7 + 9$
and $\sO_{13}, \sO_{13} + 6,$
where
\begin{gather*}
\sO_3 = \{(2, 4, 8, 8),(2, 8, 8, 8),(8, 8, 8, 10)\} \\
\sO_4 = \{(1, 2, 2, 11),(2, 2, 5, 11),(2, 2, 5, 7),(2, 5, 10, 11)\}\\
\sO_7  = \{(0, 0, 0, 10),(0, 0, 0, 2),(0, 0, 2, 4),
(0, 0, 4, 6),(0, 0, 6, 8),(0, 0, 8, 10),(0, 4, 6, 8)\}\\
\sO_{13} = \{(0, 0, 1, 11),(0, 0, 1, 3),(0, 0, 3, 5),
(0, 0, 5, 7),(0, 0, 7, 9),
(0, 0, 9, 11),(0, 1, 3, 8),\\
(0, 3, 4, 5),(0, 3, 4, 9),(0, 3, 8, 9),(0, 4, 9, 11),
(0, 7, 8, 9),(3, 4, 8, 9)\}.
\end{gather*}
It is easy to check that each of these orbits omits at least
three residue classes modulo $12$; for example
$\sO_3$ omits the eight classes $0, 1, 3, 5, 6, 7, 9, $ and $11 \pmod{12}$,
so the theorem follows from the claim.

To prove the claim, we first calculate there are
278 solutions to the hyperbolic 
Descartes equation mod $12$,
and then check how many of these lift to global integer
solutions of the hyperbolic Descartes equation.  
Of these solutions, 184 can be eliminated by consideration 
modulo 24 and 48 as in the spherical case.
However, 13 quadruples not on the above list work modulo 48.  
These 13 quadruples lie in three orbits, namely
\begin{gather}
\sO_3 + 6 = \{(2,2,2,4),(2,2,2,8),(2,2,8,10)\},\nonumber \\
\sO_3 + 8 = \{(2,4,10,10),(4,10,10,10),(8,10,10,10)\},\label{eq5.3}\\
\sO_7 + 6 = \{(0,2,6,10),(0,2,6,6),(0,6,6,10),(2,4,6,6),(4,6,6,6),
(6,6,6,8),(6,6,8,10)\}.\nonumber
\end{gather}
To eliminate these quadruples, 
we reduce the hyperbolic Descartes equation modulo 96, which yields
\begin{multline}\label{mod96}
2((12a + w)^2 + (12b + x)^2 + (12c + y)^2 + (12d + z)^2) \equiv\\
         (12(a + b + c + d) + (w + x + y + z))^2 + 4 \pmod{96}.
\end{multline}
After substituting  $(w,x,y,z) = (2,2,2,4)$, $(2,4,10,10)$, or 
$(0,2,6,10)$ into \eqn{mod96} 
we get
\beqs~ \label{mod96.2}
48\alpha^2 + 48\alpha + 48 \equiv 0 \pmod {96}
\eeqs
where $\alpha = a + b + c + d$.  This equation has 
no solutions as $\alpha ^2 + \alpha + 1 \equiv 1 \pmod 2$ 
for any $\alpha$.
We have shown that one quadruple from each of the orbits \eqn{eq5.3} 
is impossible, thus ruling out the entire orbit.

Finally, a quick search shows that the other 15 orbits all occur 
as reductions of integer hyperbolic Apollonian quadruples. ~~~$\Box$
\end{pf}

\begin{table}[htbp]
\begin{center}
$n \equiv 0 \pmod{24}$
\begin{tabular}{|cccccccccc|}
\hline 
  24&  48&  144&  168&  264&  384&  504&  528&  720&  792\\
  864&  960&  984&  1008&  1104&  1344&  1392&  1632&  1728&  1824\\
  2208&  2232&  2352&  2448&  2592&  2904&  3144&  3192&  3384&  3744\\
  3984&  4032&  4104&  4248&  4464&  4584&  4944&  5352&  5424&  5664\\
  5784&  6048&  6384&  6624&  7272&  7344&  7464&  7776&  8064&  8160\\
  8664&  8712&  8808&  8904&  9264&  9312&  9984&  10032&  10224&  10248\\
  10368&  10752&  11064&  11184&  11304&  11424&  11592&  11952&  12648&  12864\\
  13272&  13368&  13584&  13704&  13824&  14064&  14160&  14304&  15144&  15552\\
  15624&  15984&  16704&  16872&  17712&  17784&  18192&  18264&  18768&  19224\\
  19704&  19872&  19944&  20304&  20424&  20664&  22104&  22632&  23304&  23568\\
  24744&  24816&  24984&  25248&  25944&  26544&  26904&  28008&  28224&  28392\\
  29376&  29784&  29976&  30072&  30744&  32544&  32904&  33552&  33744&  33888\\
  33936&  33984&  36024&  36672&  36984&  37104&  37464&  37632&  38424&  39024\\
  39096&  39624&  40104&  40152&  41784&  41904&  41952&  42816&  43272&  43584\\
  43848&  43896&  43944&  44856&  45144&  45912&  46752&  46824&  46872&  47184\\
  47784&  49704&  49992&  50160&  50448&  50928&  51312&  51504&  51552&  52584\\
  53424&  54816&  54864&  55464&  55728&  55824&  56736&  56856&  58344&  58824\\
  60528&  60744&  61872&  62424&  63120&  64224&  64344&  64824&  64944&  65544\\
  65664&  66144&  66744&  66984&  68304&  68472&  70344&  71712&  71856&  72744\\
  73944&  76560&  76584&  77592&  77664&  78384&  79584&  80664&  81264&  81984\\
  83784&  84384&  84624&  85536&  87744&  88152&  88344&  92712&  93024&  93672\\
  93864&  94584&  95184&  96144&  96216&  97848&  99264&  99792&  99864&  99984\\\hline
\end{tabular}

\bigskip

$n \equiv 1 \pmod{24}$

\begin{tabular}{|cccccccccc|}
\hline 
 49&  217&  529&  553&  889&  1441&  1897&  2737&  3073&  3337\\
  4009&  4417&  4609&  6049&  7273&  8449&  9289&  9889&  10609&  10921\\
  11017&  11257&  11809&  11929&  12649&  14449&  15289&  18529&  18601&  20209\\
  23017&  24577&  27769&  37009&  39577&  43489&  45649&  46129&  47449&  51049\\
  52369&  54313&  54529&  55369&  56809&  66289&  66889&  73249&  80569&  83329\\
  91129&  91729&  92329&&&&&&&\\\hline
\end{tabular}

\bigskip 

$n \equiv 2 \pmod{24}$
\begin{tabular}{|cccccccccc|}
\hline 
434&  506&  1034&  3626&  7706&  8786&  12674&  19634&  25154&  27554\\
  28034&  29714&  41354&  83426&  97850&&&&&\\\hline
\end{tabular}

\bigskip 

$n \equiv 5 \pmod{24}$
\begin{tabular}{|cccccccccc|}
\hline 
 77&  749&  869&  893&  1301&  2189&  2261&  2429&  2573&  3317\\
  3509&  5789&  6077&  9437&  16469&  17789&  19589&  22589&  22709&  44549\\
  56909&  72869&&&&&&&&\\\hline

\end{tabular}
\caption{Missing curvatures below $10^{5}$ in the spherical $(0,1,1, 2)$ packing}
\label{tab:2}
 \end{center}
 \end{table}

In \cite{GLMWY21}, it is conjectured that in any 
integral Euclidean Apollonian packing, all
sufficiently large integers not ruled out
by congruence conditions modulo 24 are represented.
This was suggested
by numerical data for a few specific packings. 
Table~\ref{tab:2} presents
 integers below 100000 missing from  the 
spherical Apollonian packing $(0,1,1,2)$ grouped by residue
classes modulo $24$; we give data for 
residue classes $0, 1, 2, 5 \pmod{24}$ as
representative. It is not completely clear from this data whether
one should believe
that only a finite number of integers will be missed in each class.
For certain residue classes modulo $24$ the exceptions seem
to be rapidly thinning out. However for the class $0~ (\bmod~24)$
(and also  $9$, $18$ and $21~ (\bmod~24)$) the data seems
equivocal.

It seems very likely that there are no congruence restrictions on integers
in spherical and hyperbolic circle packings for any modulus 
prime to $6$. 
In \cite[Theorem 6.2]{GLMWY21}, it is shown that
in any integral Euclidean Apollonian packing,
the ``curvatures'' of the circles include representatives of
every residue class modulo $m$ for any modulus $m$ that is
relatively prime to 30.  We expect a similar  result to hold for
spherical and hyperbolic packings as well.

{\tt
\begin{tabular}{ll}
email: &eriksson@math.berkeley.edu \\
& lagarias@umich.edu \\
\end{tabular}
 }
\end{document}